# ASYMPTOTIC EXPANSIONS FOR SUMS OF BLOCK-VARIABLES UNDER WEAK DEPENDENCE[1]

BY S. N. LAHIRI

*Iowa State University*

Let $\{X_i\}_{i=-\infty}^{\infty}$ be a sequence of random vectors and $Y_{in} = f_{in}(\mathcal{X}_{i,\ell})$ be zero mean block-variables where $\mathcal{X}_{i,\ell} = (X_i, \ldots, X_{i+\ell-1}), i \geq 1$, are overlapping blocks of length $\ell$ and where $f_{in}$ are Borel measurable functions. This paper establishes valid joint asymptotic expansions of general orders for the joint distribution of the sums $\sum_{i=1}^{n} X_i$ and $\sum_{i=1}^{n} Y_{in}$ under weak dependence conditions on the sequence $\{X_i\}_{i=-\infty}^{\infty}$ when the block length $\ell$ grows to infinity. In contrast to the classical Edgeworth expansion results where the terms in the expansions are given by powers of $n^{-1/2}$, the expansions derived here are mixtures of two series, one in powers of $n^{-1/2}$ and the other in powers of $[\frac{n}{\ell}]^{-1/2}$. Applications of the main results to (i) expansions for Studentized statistics of time series data and (ii) second order correctness of the blocks of blocks bootstrap method are given.

**1. Introduction.** Let $\{X_i\}_{i \in \mathbb{Z}}$ be a sequence of $\mathbb{R}^{d_0}$-valued $(d_0 \in \mathbb{N})$ random vectors such that $EX_i = 0$ for all $i \in \mathbb{Z}$, where $\mathbb{Z} = \{0, \pm 1, \pm 2, \ldots\}$ and $\mathbb{N} = \{1, 2, \ldots\}$. The sequence $\{X_i\}_{i \in \mathbb{Z}}$ need not be stationary. Let $\chi_{i,\ell} \equiv (X_i, \ldots, X_{i+\ell-1})', i, \ell \in \mathbb{N}$, denote (overlapping) blocks of length $\ell$ for some given integer $\ell \equiv \ell_n \in [1, n]$ and let

(1.1) $$Y_{in} = f_{in}(\chi_{i,\ell}), \qquad i, n \in \mathbb{N},$$

denote the block-variables, where $f_{in} : \mathbb{R}^{d_0 \ell} \to \mathbb{R}^{d_1}, d_1 \in \mathbb{N}$, are Borel-measurable functions such that $EY_{in} = 0$ for all $i, n \in \mathbb{N}$. Let $b \equiv \lceil n/\ell \rceil$, where for any $x \in \mathbb{R}$, $\lceil x \rceil$ denotes the smallest integer not less than $x$. The main results of

Received March 2005; revised June 2006.
[1]Supported in part by NSF Grants DMS-00-72571 and DMS-03-06574.

the paper give asymptotic expansions for the scaled sums

$$(1.2) \qquad S_n = \left(\frac{1}{\sqrt{n}}\sum_{i=1}^{n} X_i', \frac{1}{\sqrt{n\ell}}\sum_{i=1}^{n} Y_{in}'\right)' \equiv (S_{1n}', S_{2n}')', \qquad n \geq 1,$$

under some weak dependence conditions on $\{X_i\}_{i\in\mathbb{Z}}$ when $\ell_n$ becomes unbounded as $n \to \infty$. Here and in the following, $d \equiv d_0 + d_1$ denotes the dimension of $S_n$, and $A'$ denotes the transpose of a matrix $A$. The block-variables $Y_{in}$ serve as basic building blocks for many important statistical methods used in the analysis of weakly dependent time series data. Some important examples are given below.

EXAMPLE 1.1 (*Spectral density estimation*). Suppose that $\{X_i\}_{i\in\mathbb{Z}}$ is a second-order stationary process with values in $\mathbb{R}$ such that $\sum_{k=0}^{\infty} |\gamma(k)| < \infty$, where $\gamma(k) = \text{Cov}(X_1, X_{1+k}), k \geq 0$. Then $\{X_i\}_{i\in\mathbb{Z}}$ has a spectral density $f : (-\pi, \pi) \to [0, \infty)$. A general class of nonparametric estimators of $f$ (cf. Priestley [26]) is given by

$$(1.3) \qquad \hat{f}_n(\lambda) = (2\pi)^{-1} \sum_{k=0}^{\ell} \omega_{kn} \hat{\gamma}_n(k) \cos(k\lambda), \qquad \lambda \in (-\pi, \pi),$$

where $\omega_{kn} \in \mathbb{R}$ are nonrandom weights and where for $0 \leq k \leq n-1$, $\hat{\gamma}_n(k) \equiv n^{-1} \sum_{i=1}^{n-k} X_i X_{i+k} - \bar{X}_n^2$ is a version of the sample lag-$k$ autocovariance. The estimator $\hat{f}_n$ plays a fundamental role in the frequency domain analysis of time series data; see, for example, Priestley [26]. Note that $\hat{f}_n$ of (1.3) can be expressed as a function of the sum $S_n$ of (1.2) where the block variables $Y_{in}$ are given by $Y_{in} = X_i[\sum_{k=0}^{\ell} \tilde{\omega}_{ikn} X_{i+k} \cos(k\lambda)], 1 \leq i \leq n$, for some suitable constants $\tilde{\omega}_{ikn} \in \mathbb{R}$, depending on the $\omega_{kn}$'s.

EXAMPLE 1.2 (*Block bootstrap methods*). Let $\hat{\theta}_n = t_n(X_1, \ldots, X_n)$ be an estimator of a parameter of interest $\theta \in \mathbb{R}$ where $t_n$ is a symmetric function of its arguments. For estimating the distribution of $\hat{\theta}_n$, Künsch [16] and Liu and Singh [23] proposed the moving block bootstrap (MBB) method. We now briefly describe the MBB for later reference. Let $\ell \equiv \ell_n \in (1, n)$ be a given integer such that (for simplicity) $n/\ell \in \mathbb{N}$, and let $\chi_{1,\ell}^*, \ldots, \chi_{b,\ell}^*$ be selected at random, with replacement from the "observed" blocks $\{\chi_{1,\ell}, \ldots, \chi_{N,\ell}\}$, where $b \equiv \lceil n/\ell \rceil = n/\ell$ and $N = n - \ell + 1$. Let $\theta_n^* \equiv t_n(X_1^*, \ldots, X_n^*)$ denote the MBB version of $\hat{\theta}_n$, where $X_1^*, \ldots, X_n^*$ are elements of the resampled blocks $\chi_{i,\ell}^*, 1 \leq i \leq b$. Then the MBB estimator of the distribution function $G_n(x) \equiv P(\hat{\theta}_n \leq x)$ of $\hat{\theta}_n$ is given by

$$\hat{G}_n(x) = P(\theta_n^* \leq x | X_1, \ldots, X_n),$$

the conditional distribution function of $\theta_n^*$, given $X_1, \ldots, X_2$, and the MBB estimator of a functional of $G_n$ is given by "plugging in" $\hat{G}_n$ for $G_n$. Since $\chi_{1,\ell}^*, \ldots, \chi_{b,\ell}^*$ are conditionally independent and identically distributed (i.i.d.) with $\Pr(\chi_{1,\ell}^* = \chi_{i,\ell}) = N^{-1}$, $1 \leq i \leq N$, it follows that $\hat{G}_n(x)$ and its functionals can be represented as functions of the block-variables $Y_{in} = f_{in}(\chi_{i,\ell}), 1 \leq i \leq N$, for suitable functions $f_{in}$ (with $Y_{in} = 0$ for $N+1 \leq i \leq n$). For example, if $d_0 = 1$ and $\hat{\theta}_n = n^{-1} \sum_{i=1}^n X_i$ is the sample mean, then it is easy to check that the MBB estimator $\hat{\sigma}_n^2$ of the variance of $\hat{\theta}_n$ is given by

$$(1.4) \qquad \hat{\sigma}_n^2 = n^{-1} \left[ N^{-1} \sum_{i=1}^N U_{1i}^2 - \left( N^{-1} \sum_{i=1}^N U_{1i} \right)^2 \right],$$

where $U_{1i} \equiv U_{1in} = (X_i + \cdots + X_{i+\ell-1})/\sqrt{\ell}$ is the scaled sum of the $i$th block $\chi_{i,\ell}$.

The subsampling method of Politis and Romano [25] and Hall and Jing [13], and the block empirical likelihood method of Kitamura [15] are other important examples of resampling methods that naturally lead to sums of block-variables.

EXAMPLE 1.3 (*Studentized statistics*). Suppose that $d_0 = 1$ and $\{X_i\}_{i \in \mathbb{Z}}$ is stationary. Let $\hat{\theta}_n = H(\bar{X}_n)$ be an estimator of $\theta = H(EX_1)$, where $H : \mathbb{R} \to \mathbb{R}$ is a smooth function. This is a version of the "smooth function" model (cf. Hall [11]) that covers many commonly used estimators. For constructing confidence intervals for $\theta$, one considers the approximate pivotal statistic

$$(1.5) \qquad T_n = \sqrt{n}(\hat{\theta}_n - \theta)/\tilde{\tau}_n,$$

where $\tilde{\tau}_n^2$ is an estimator of the asymptotic variance of $\sqrt{n}(\hat{\theta}_n - \theta)$. For example, we may use $\tilde{\tau}_n^2 = h(\bar{X}_n)\tilde{\sigma}_n^2$ where $h(\cdot)$ denotes the derivative of $H(\cdot)$ and $\tilde{\sigma}_n^2$ is either $2\pi \hat{f}_n(0)$ of Example 1.1 or the scaled bootstrap estimator $n\hat{\sigma}_n^2$ of Example 1.2 above. In both cases, $T_n$ is a function of the sum $S_n$.

The examples above show that the scaled sum $S_n$ of (1.2) plays a fundamental role in statistical inference for weakly dependent processes. In a seminal paper, Götze and Hipp [8] derived asymptotic expansions for the scaled sum $S_{1n}$ of weakly dependent random vectors under an exponential mixing condition. This paper builds upon the work of Götze and Hipp [8] and proves asymptotic expansions for the augmented sum $S_n$, under a similar general framework. The proofs of the main results are based on some extensions and refinements of the techniques developed by Götze and Hipp [8] and Lahiri [17].

To highlight some of the major differences between the present problem and the case of the regular sum $\sum_{i=1}^n X_i$ treated by Götze and Hipp [8] and others, note that even when the $X_i$'s are strongly mixing at an exponential rate (e.g., as in [8]), the block-variables $\{Y_{in}\}_{i=1}^n$, being defined on *overlapping* blocks of length $\ell$, in general have a strong mixing coefficient equal to *one* for all lags of order $\leq (\ell - 1)$. Since $\ell \to \infty$ with $n$ [in this paper, $\ell$ could grow as fast as $O(n^{1-\kappa})$ for a given $\kappa \in (0, 1)$], this leads to a very strong form of dependence among an *unbounded* number of neighboring block variables $Y_{in}$'s, thereby destroying the weak dependence structure of the original sequence $\{X_i\}_{i \in \mathbb{Z}}$. The "local strong dependence" of the $Y_{in}$'s has a nontrivial effect on the accuracy of approximation and on the growth rate of the variance of $\sum_{i=1}^{n} Y_{in}$. Indeed, the sum of $Y_{in}$'s over a block of size $\ell$ is of the order $O_P(\ell)$ compared to the order $O_P(\ell^{1/2})$ for weakly dependent variables and hence, the componentwise variance terms of $\sum_{i=1}^{n} Y_{in}$ typically grow at the rate $O(n\ell)$. This leads to the normalizing constant $(n\ell)^{-1/2}$ for the sum $\sum_{i=1}^{n} Y_{in}$ in (1.2). Intuitively, the "local strong dependence" of the $Y_{in}$'s makes the sum $\sum_{i=1}^{n} Y_{in}$ essentially behave like a sum of $O(n/\ell)$-many "approximately independent" variables. As a result, the accuracy of an $(s-2)$-order asymptotic expansion for $S_n$ is only $o([n/\ell]^{-(s-2)/2})$, which should be compared to the order $o(n^{-(s-2)/2})$ for $S_{1n}$. Further, in contrast to the case of $S_{1n}$, *asymptotic expansions for $S_n$ of* (1.2) *are now given by a mixture of two series of terms, one in powers of* $n^{-1/2}$ (*corresponding to* $\frac{1}{\sqrt{n}} \sum_{i=1}^{n} X_i$) *and the other in powers of* $[n/\ell]^{-1/2}$ (*corresponding to* $\frac{1}{\sqrt{n\ell}} \sum_{i=1}^{n} Y_{in}$).

The main results of the paper give an $(s-2)$th order expansion for $Ef(S_n)$ for Borel-measurable functions $f : \mathbb{R}^d \to \mathbb{R}$ for any integer $s \geq 3$. For *smooth* functions $f$, expansions for $Ef(S_n)$ are established *without* requiring any type of Cramér's condition. For i.i.d. random vectors, such results were first proved by Götze and Hipp [7] for the regular sum $\sum_{i=1}^{n} X_i$. Further, moderate deviation bounds for $S_n$ are also proved in Section 2 below. In each of these results, the block length variable $\ell$ is allowed to go to infinity at a rate $O(n^{1-\kappa})$ for an arbitrarily small $\kappa > 0$.

Three important applications of the main results are considered in Section 3. The first result establishes moderate deviation bounds for the MBB moments. These bounds are useful for studying accuracy of MBB variance and distribution function estimators. The second result gives a *general order* Edgeworth expansion (EE) for a version of the Studentized sample mean, where the conditional Cramér condition is verified *explicitly*. It may be noted that the standard EE theory based on sums of *finite*-dimensional random vectors has severe limitations in this problem, as the Studentizing factor in the dependent case is no longer a smooth function of finitely many sample means. A second notable feature of this EE result is that the $(s-2)$-order $(s \geq 3)$ expansion for the Studentized sample mean is proved here *solely* under a conditional Cramér condition on the variables $X_i$'s. This is in sharp contrast to the i.i.d. case, where some additional conditions, like a Cramér condition on the joint distribution of $(X_1, X_1^2)'$, are required (cf. Bhattacharya and Ghosh [2]). Intuitively, here the block variables become subjected to a central limit theorem (CLT) effect due to block averaging, and this entails the required conditional Cramér condition for the joint distribution of the linear and quadratic parts.

The third application of the main results of Section 2 is to studying higher-order properties of the blocks of blocks bootstrap (BOBB) method of Politis and Romano [24]. It is shown that under some regularity conditions, the BOBB approximation to a class of Studentized statistics is second-order correct (s.o.c.), that is, it is more accurate than the limiting normal distribution. This result is particularly useful for constructing s.o.c. confidence intervals for infinite-dimensional parameters of the underlying process, such as the spectral density.

The rest of the paper is organized as follows. Section 1 concludes with a brief literature review. Section 2 gives the main results on expansions for $S_n$, while Section 3 gives results on the MBB moments, the Studentized sample mean, and the BOBB method. Proofs of all the results are given in Section 4.

There is a vast literature on asymptotic expansions for $S_{1n}$ and for statistics that are smooth functions of $S_{1n}$. A detailed account of the theory for sums of independent random vectors is given in Bhattacharya and Ranga Rao [3]. The theory under the "smooth function model" is treated by Bhattacharya and Ghosh [2], Skovgaard [27] and Hall [11], among others. For sums of *weakly dependent* random vectors, Götze and Hipp [8] obtained expansions under a very flexible framework. Applicability of [8] results in different time series models has been verified in [9]. Lahiri [17] relaxes the moment condition used by Götze and Hipp [8] and settles a conjecture of [8] on the validity of expansions for expectations of smooth functions of $S_{1n}$. Expansions for $S_{1n}$ under polynomial mixing rates have been given by Lahiri [18]. EEs for Studentized statistics under weak dependence are given by Götze and Künsch [10] and Lahiri [19] to the second order for general weakly dependent processes, and by Velasco and Robinson [30] to higher orders for Gaussian processes. EEs of a general order without the Gaussian assumption have recently been established in Lahiri [21].

**2. Main results.** Let $\{X_i\}_{i\in\mathbb{Z}}$ be a sequence of (possibly nonstationary) zero mean $\mathbb{R}^{d_0}$-valued random vectors (as in Section 1) and let $Y_{in}, i \geq 1, n \geq 1$, be as defined in (1.1). In this section, we establish asymptotic expansions for the scaled sum $S_n$ of (1.2) under a framework similar to [8]. Suppose that the $X_i$'s are defined on a probability space $(\Omega, \mathcal{F}, P)$ and that $\{\mathcal{D}_j\}_{j=-\infty}^{\infty}$ is a given collection of sub-$\sigma$-fields of $\mathcal{F}$. Let $\mathcal{D}_p^q \equiv \sigma\langle\{\mathcal{D}_j : j \in \mathbb{Z}, p \leq j \leq q\}\rangle, -\infty \leq p < q \leq \infty$. For $k = 1, \ldots, b$, define

$$\bar{X}_{kn} = \ell^{-1} \sum_{i=(k-1)\ell+1}^{k\ell \wedge n} X_i, \qquad \bar{Y}_{kn} = \ell^{-1} \sum_{i=(k-1)\ell+1}^{k\ell \wedge n} Y_{in},$$

(2.1)
$$W_{kn} = (\sqrt{\ell}\bar{X}'_{kn}; \bar{Y}'_{kn})',$$

where $b \equiv \lceil n/\ell \rceil$ and $x \wedge y = \min\{x, y\}$, $x, y \in \mathbb{R}$. Let $\|x\| = (x_1^2 + \cdots + x_k^2)^{1/2}$ and $|x| = |x_1| + \cdots + |x_k|$, $x = (x_1, \ldots, x_k)' \in \mathbb{R}^k$ and let $|B|$ denote the size of a set $B$. Let $\tilde{b} = n/\ell$. For notational simplicity, we drop the subscript $n$ in $\ell, b, \tilde{b}$. Unless otherwise stated, the limits in the order symbols are taken by letting $n \to \infty$.

CONDITIONS.

C.1. There exists $\kappa \in (0, 1)$ such that $\kappa \log n < \ell < \kappa^{-1} n^{1-\kappa}$ for all $n > \kappa^{-1}$.

C.2. (i) There exist $\rho \in (0, \infty)$ and $s \in \{3, 4, \ldots\}$ such that

(2.2) $\max\{Eh_s(\|\bar{Y}_{kn}\|), Eh_s(\|\sqrt{\ell}\bar{X}_{kn}\|)\} \leq \rho$ for all $1 \leq k \leq b, n \geq 1$,

where $h_s(u) \equiv u^s[\log(1+u)]^{\alpha(s)}, u \geq 0$, and $\alpha(s) = 2s^2$.

(ii) For all $i \geq 1, n \geq 1$, $EX_i = 0 = EY_{in}$, and

(2.3) $\lim_{n\to\infty} \text{Cov}(S_n) = \Xi_\infty$ exists and is nonsingular.

Further, there exists a $\kappa \in (0, 1)$ such that $\inf_{\|t\|=1} t' \text{Cov}(\sum_{j=j_0+1}^{j_0+m} W_{jn})t > \kappa m$ for all integers $j_0, m, n$ with $0 \leq j_0 \leq b - m$, $\kappa^{-1} \leq m \leq \sqrt{b}$ and $n \geq \kappa^{-1}$.

C.3. There exists $\kappa \in (0, 1)$ such that for $m, n > \kappa^{-1}$ and for $j \geq 1, 1 \leq i \leq n$, there exist a $\mathcal{D}_{j-m}^{j+m}$-measurable $X_{j,m}^\dagger$ and a $\mathcal{D}_{i-m}^{i+m+\ell}$-measurable $Y_{in,m}^\dagger$ such that

(2.4) $\max\{E\|X_j - X_{j,m}^\dagger\|, \ell^{-1}E\|Y_{in} - Y_{in,m}^\dagger\|\} \leq \kappa^{-1}\exp(-\kappa m)$.

C.4. There exists a constant $\kappa \in (0, 1)$ such that for all $m \in \mathbb{N}$,

$\sup\{|P(A \cap B) - P(A)P(B)| : A \in \mathcal{D}_{-\infty}^i, B \in \mathcal{D}_{i+m}^\infty, i \in \mathbb{Z}\} \leq \kappa^{-1}\exp(-\kappa m)$.

C.5. There exists a constant $\kappa \in (0, 1)$ such that for all $i, j, k, r, m = 1, 2, \ldots$ and $A \in \mathcal{D}_i^j$ with $i < k < r < j$ and $m > \kappa^{-1}$,

$$E|P(A|\mathcal{D}_j : j \notin [k, r]) - P(A|\mathcal{D}_j : j \in [i-m, k] \cup (r, j+m])|$$
$$\leq \kappa^{-1}\exp(-\kappa m).$$

C.6. There exist $a \in (0, \infty), \kappa \in (0, 1)$ and sequences $\{m_n\} \subset \mathbb{N}$ and $\{d_n\} \subset [1, \infty)$ with $m_n^{-1} + m_n b^{-1/2} = o(1), d_n = O(\ell + b^a)$ and $d_n^2 m_n = O(b^{1-\kappa})$ such that

$$\max_{j_0 \in J_n} \sup_{t \in A_n} E\left|E\left\{\exp\left(\iota t' \sum_{j=j_0-m_n}^{j_0+m_n} W_{jn}\right)\bigg|\tilde{\mathcal{D}}_{j_0}\right\}\right| \leq 1 - \kappa \quad \text{for all } n \geq \kappa^{-1},$$
(2.5)

where $J_n = \{m_n + 1, \ldots, b - m_n + 1\}, A_n = \{t \in \mathbb{R}^d : \kappa d_n \leq \|t\| \leq [b^a + \ell]^{1+\kappa}\}$, and $\tilde{\mathcal{D}}_{j_0} = \sigma\langle\{\mathcal{D}_j : j \in \mathbb{Z}, j \notin [(j_0 - \lfloor\frac{m_n}{2}\rfloor)\ell + 1, (j_0 + \lfloor\frac{m_n}{2}\rfloor + 1)\ell]\}\rangle$.

Now we comment on the conditions. C.1 is a growth condition on the block length $\ell$ and covers optimal block sizes in most applications. For example, for an exponentially strongly mixing process, the optimal block size for the spectral density estimation problem (cf. Example 1.1) is $O(\log n)$, while for the block bootstrap estimation of variance and distribution functions (cf. Example 1.2), the optimal block lengths are of the order $n^{1/k}$ for $k = 3, 4, 5$ (cf. Hall, Horowitz and Jing [12]), all of which are covered by C.1. Next consider C.2–C.6. As in [8], here we formulate the conditions in terms of the auxiliary $\sigma$-fields $\mathcal{D}_j$ in order to allow for greater generality. Condition C.2(i) is a moment condition on $\sqrt{\ell}\bar{X}_{kn}$ and $\bar{Y}_{kn}$ that is optimal up to the logarithmic factor. Since $h_s(\cdot)$ is convex and nondecreasing, $Eh_s(\|\bar{Y}_{kn}\|) \leq \ell^{-1}\sum_{i=(k-1)\ell+1}^{k\ell} Eh_s(\|Y_{in}\|)$, and a sufficient condition for (2.2) is $\max\{Eh_s(\|Y_{in}\|): 1 \leq i \leq n, n \geq 1\} \leq \rho$. Equation (2.3) ensures a nondegenerate normal limit for $S_n$ with the given scaling constants. If the process $\{X_i\}_{i=-\infty}^{\infty}$ is second-order stationary and if, for each $n \geq 1$, the collection $\{Y_{in}\}_{i=1}^{n}$ is also second-order stationary, then the last part of C.2(ii) follows from (2.3). Condition C.3 connects the variables $X_i$ and $Y_{in}$ to the strong-mixing property C.4 of the auxiliary $\sigma$-fields $\mathcal{D}_j$. For suitably well-behaved functions $f_{in}$, the bound on $E\|Y_{in} - Y_{in,m}^{\dagger}\|$ follows from the bounds on the $E\|X_i - X_{i,m}^{\dagger}\|$'s. C.5 is an approximate Markov property.

Finally, consider the conditional Cramér condition C.6, which is weaker than a direct analog of the conditional Cramér condition of [8] for the summands in $S_n$ on two counts. First, the range of values of $t$ in (2.5) has a lower bound, namely, $\kappa d_n$, that may tend to infinity at a suitable rate. By setting $d_n \equiv 1$ for all $n \geq 1$, one gets the more standard version of the Cramér condition, as formulated by Götze and Hipp [8]. Second, the conditioning $\sigma$-field $\tilde{\mathcal{D}}_{j_0}$ in (2.5) is only a *sub*-$\sigma$-field of the standard choice $\mathcal{D}_{-j_0} \equiv \sigma\langle\{\mathcal{D}_j : j \in \mathbb{Z}, j \neq j_0\}\rangle$ (cf. (2.5) of [8]). Indeed, if $\mathcal{D}_{j_0}^* \subset \mathcal{F}$ is a $\sigma$-field of $\mathcal{F}$ containing $\tilde{\mathcal{D}}_{j_0}$, then one can show (cf. Lahiri [22]) that a sufficient condition for (2.5) is given by

$$(2.6) \qquad \max_{j_0 \in J_n} \sup_{t \in A_n} E\left|E\left\{\exp\left(\iota t' \sum_{j=j_0-m_n}^{j_0+m_n} W_{jn}\right)\Big|\mathcal{D}_{j_0}^*\right\}\right| \leq 1 - \kappa.$$

Indeed, these two refinements allow us to establish valid EEs for the *Studentized* sample mean solely under Götze and Hipp's [8] Cramér condition on the $X_i$'s (cf. Section 3.2).

Next define the functions $P_{r,n}(t)$ for $t \in \mathbb{R}^d$ by the identity (in $u \in \mathbb{R}$)

$$(2.7) \qquad \exp\left(\sum_{r=3}^{s}(r!)^{-1}u^{r-2}\tilde{b}^{(r-2)/2}\chi_{r,n}(t)\right) = 1 + \sum_{r=1}^{\infty} u^r P_{r,n}(t),$$

where $\chi_{r,n}(t)$ denotes the $r$th cumulant of $t'S_n$ [multiplied by $(\iota)^r$], defined by

$$\chi_{r,n}(t) = \frac{d^r}{du^r} \log E \exp(\iota u t' S_n)\bigg|_{u=0}. \tag{2.8}$$

Under C.1–C.5, the $P_{r,n}(t)$'s are bounded for each $t$. Define the $(s-2)$th order EE $\Psi_{s,n}$ of $S_n$ through its Fourier transform $\hat{\Psi}_{s,n}(t) \equiv \int e^{\iota t'x} d\Psi_{s,n}(x), t \in \mathbb{R}^d$, by

$$\hat{\Psi}_{s,n}(t) = \exp(\chi_{2,n}(t)/2)\left[1 + \sum_{r=1}^{s-2} \tilde{b}^{-r/2} P_{r,n}(t)\right], \qquad t \in \mathbb{R}^d. \tag{2.9}$$

Next, let $s_0 = 2\lfloor s/2 \rfloor$ and for any positive definite matrix $A$ of order $k \in \mathbb{N}$, let $\Phi_A$ and $\Phi(\cdot; A)$ both denote the normal distribution on $\mathbb{R}^k$ with mean zero and covariance matrix $A$. Then we have the following result on expansions for $Ef(S_n)$.

THEOREM 2.1. *Assume that conditions* C.1–C.6 *hold for some* $a > (s-2)/2$ *in* C.6. *Let* $f : \mathbb{R}^d \to \mathbb{R}$ *be a Borel measurable function with* $M_f \equiv \sup\{(1+\|x\|^{s_0})^{-1}|f(x)| : x \in \mathbb{R}^d\} < \infty$. *Then there exist constants* $C_1 = C_1(a)$, $C_2 \in (0, \infty)$ *(neither depending on* $f$*) such that for all* $n > C_2$,

$$\left|Ef(S_n) - \int f \, d\Psi_{s,n}\right| \leq C_1 \omega(\tilde{f}; b^{-a}) + C_2 M_f b^{-(s-2)/2} (\log n)^{-2}, \tag{2.10}$$

*where* $\tilde{f}(x) = f(x)/(1 + \|x\|^{s_0}), x \in \mathbb{R}^d$, *and* $\omega(\tilde{f}; \varepsilon) = \int \sup\{|\tilde{f}(x+y) - \tilde{f}(x)| : \|y\| \leq \varepsilon\} \Phi(dx; \Xi_\infty), \varepsilon > 0$.

A direct consequence of Theorem 2.1 is that $\sup_{A \in \mathcal{C}} |P(S_n \in A) - \Psi_{s,n}(A)| = O(b^{-(s-2)/2}(\log n)^{-2})$, where $\mathcal{C}$ is the collection of all measurable convex sets in $\mathbb{R}^d$.

Next consider expansions of $Ef(S_n)$ when $f$ is smooth. When the $X_i$'s are i.i.d., Götze and Hipp [7] established an $(s-2)$th order expansion of $Ef(S_{1n})$ for $f \in C^{s-2}(\mathbb{R}^{d_0})$, *without* any Cramér type conditions, where $C^r(\mathbb{R}^k)$ denotes the set of all $r$-times continuously differentiable functions from $\mathbb{R}^k$ to $\mathbb{R}$, $r \geq 0, k \in \mathbb{N}$. For weakly dependent $X_i$'s, Lahiri [17] proved a similar result for $Ef(S_{1n})$ for $f \in C^{s-1}(\mathbb{R}^{d_0})$, settling a conjecture of Götze and Hipp [8]. In the same spirit, Theorem 2.2 below establishes validity of an $(s-2)$th order expansion for $Ef(S_n)$ for $f \in C^{s-1}(\mathbb{R}^d)$, without the conditional Cramér condition C.6. To state it, let $D^\alpha = \partial^{|\alpha|}/\partial t_1^{a_1} \cdots \partial t_k^{a_k}$, for $\alpha = (\alpha_1, \ldots, \alpha_k)' \in (\mathbb{Z}_+)^k$, $k \in \mathbb{N}$, where $\mathbb{Z}_+ = \{0, 1, \ldots\}$.

THEOREM 2.2. *Assume that conditions* C.1–C.5 *hold. Let* $f : \mathbb{R}^d \to \mathbb{R}$ *be a function such that* (i) $f \in C^{s-1}(\mathbb{R}^d)$ *and* (ii) *for each* $\alpha \in (\mathbb{Z}+)^d$ *with* $0 \le |\alpha| \le s-1$, there exists a $p(\alpha) \in \mathbb{Z}_+$, with $p(0) = s_0$, such that $M_{f,\alpha} \equiv \sup\{|D^\alpha f(x)|/(1+\|x\|^{p(\alpha)}) : x \in \mathbb{R}^d\} < \infty$. Then there exists a constant $C_3 \in (0, \infty)$, depending only on $s, d, \rho, \kappa$ and on $\{M_{f,\alpha} : |\alpha| \le s-1\}$, such that

$$\left| Ef(S_n) - \int f \, d\Psi_{s,n} \right| \le C_3 b^{-(s-2)/2} (\log n)^{-2} \qquad \text{for all } n \ge 2.$$

The next result is a moderate deviation inequality for $S_n$ and its moments up to order $s_0$. Let $\mathbb{1}(B)$ denote the indicator of a set $B$.

THEOREM 2.3. *Let $\lambda_0$ denote the largest eigenvalue of $\Xi_\infty$. Then under conditions C.1–C.5, for any $\lambda > \lambda_0$, there exists a constant $C_4 \in (0, \infty)$, depending only on $s, d, \lambda, \rho, \kappa$, such that for all $n \ge 2$,*

(2.11)
$$E(1 + \|S_n\|^{s_0})\mathbb{1}(\|S_n\| > [(s-2)\lambda \log n]^{1/2})$$
$$\le C_4 b^{-(s-2)/2} (\log n)^{-2}.$$

For the regular sum $S_{1n}$, *sharper* bounds are available (cf. [8, 17]). Let $\lambda_1$ be the largest eigenvalue of $\Xi_{1,1} = \lim_{n \to \infty} \text{Cov}(S_{1n})$. Then Lahiri [17] proved that $E(1 + \|S_{1n}\|^{s_0})\mathbb{1}(\|S_{1n}\| > [(s-2)\lambda \log n]^{1/2}) = o(n^{-(s-2)/2})$ for $\lambda > \lambda_1$. For the marginal distribution of $S_{2n}$, the following bound is slightly better than (2.11) (as $\lambda_2 \le \lambda_0$).

THEOREM 2.4. *Suppose that $\Xi_{2,2} \equiv \lim_{n \to \infty} \text{Cov}(S_{2n})$ exists and is nonsingular, and conditions C.1–C.5 hold for $\bar{Y}_{jn}$'s only. Then, for any $\lambda > \lambda_2$, there exists a constant $C_5 \in (0, \infty)$, depending only on $s, d, \lambda, \rho, \kappa$, such that for all $n \ge 2$,*

(2.12)
$$E(1 + \|S_{2n}\|^{s_0})\mathbb{1}(\|S_{2n}\| > [(s-2)\lambda \log n]^{1/2}) \le C_5 b^{-(s-2)/2} (\log n)^{-2},$$

*where $\lambda_2$ is the largest eigenvalue of $\Xi_{2,2}$.*

REMARK. Analogs of Theorems 2.1–2.4 are known for the sum $S_n^{[2]} = (\frac{1}{\sqrt{n}} \sum_{i=1}^n X_i', \frac{1}{\sqrt{b}} \sum_{i=1}^b (Y_{in}^{[2]})')'$, where $Y_{jn}^{[2]} \equiv Y_{(j-1)\ell+1,n}$, $j \in \mathbb{N}$, denote the block-variables based on *nonoverlapping* blocks of length $\ell$. See [22] for details.

## 3. Applications.

3.1. *Moderate deviation bounds for MBB moments.* Let $\{X_i\}_{i \in \mathbb{Z}}$ be an $\mathbb{R}^{d_0}$-valued stationary process. As described in Example 1.2, let $U_{1i} = (X_i + \cdots + X_{i+\ell-1})/\sqrt{\ell}$ denote the scaled-sum of the $X_i$'s in the ith block $\chi_{i,\ell} = (X_i, \ldots, X_{i+\ell-1})$, $i \geq 1$. For $\nu \in (\mathbb{Z}_+)^{d_0}$, we define the $\nu$th MBB moment as

$$\hat{\mu}_n(\nu) \equiv N^{-1} \sum_{j=1}^{N} U_{1j}^{\nu}, \tag{3.1}$$

where $N = n - \ell + 1$. For estimators $\hat{\theta}_n$ that can be represented as smooth functions of means, the MBB estimators of $\text{Var}(\hat{\theta}_n)$ and the distribution of Studentized $\hat{\theta}_n$ can be approximated through functions of the MBB moments. For example, the MBB variance estimator for the sample mean $\hat{\theta}_n = \bar{X}_n$ for $d_0 = 1$ is given by $\hat{\sigma}_n^2 = n^{-1}(\hat{\mu}_n(2) - [\hat{\mu}_n(1)]^2)$ [cf. (1.4)]. Similarly, the leading terms in the two-term EEs of the MBB distribution function estimators of the normalized and Studentized sample means are rational functions of $\hat{\mu}_n(\nu), 1 \leq |\nu| \leq 3$ (cf. Lahiri [19, 20], Götze and Künsch [10]). As a result, moderate deviation bounds on the MBB moments are important for investigating accuracy of such MBB estimators. The following result gives a moderate deviation bound for $\hat{\mu}_n(\nu)$.

THEOREM 3.1. *Suppose that $\nu \in \mathbb{N}^{d_0}$, $EX_1 = 0$, and that conditions C.1, C.4, C.5 hold and that condition C.3 holds for the $X_i$'s only. Also, suppose that $\lambda_3(\nu) \equiv \lim_{n \to \infty} \text{Var}(b^{1/2} \hat{\mu}_n(\nu)) \in (0, \infty)$ and that*

$$\rho \equiv E(\|X_1\|^{2q} \{\log(1 + \|X_1\|)\}^{\gamma(s)}) < \infty \tag{3.2}$$

*for some $\gamma(s) > s^2$ and $q \in \mathbb{N}$ with $2q > s|\nu|$. Then, for any $\lambda > \lambda_3(\nu)$, there exists a constant $C_6 = C_6(\nu, d, s, \kappa, \lambda) \in (0, \infty)$ such that for all $n \geq 2$*

$$\begin{aligned}
&E[\{1 + |\sqrt{b}(\hat{\mu}_n(\nu) - EU_{11}^{\nu})|^{s_0}\} \\
&\quad \times \mathbb{1}(|\sqrt{b}(\hat{\mu}_n(\nu) - EU_{11}^{\nu})| > [(s-2)\lambda \log n]^{1/2})] \\
&\leq C_6 b^{-(s-2)/2} (\log n)^{-2}.
\end{aligned} \tag{3.3}$$

For an exact expression for $\lambda_3(\nu)$, see [22]. The following result is a simple consequence of Theorem 3.1 and serves to illustrate how the bounds on the $\hat{\mu}_n(\nu)$'s can be used for deriving similar bounds on the MBB estimators of population parameters, such as on the MBB estimator $\hat{\sigma}_n^2 = n^{-1}(\hat{\mu}_n(2) - [\hat{\mu}_n(1)]^2)$ of $\text{Var}(\bar{X}_n)$.

COROLLARY 3.1. *Suppose that $d_0 = 1$ and that the conditions of Theorem 3.1 hold with $\nu = 2$. Let $\hat{\Delta}_n = n[\hat{\sigma}_n^2 - E\hat{\sigma}_n^2]$. Then for any $\lambda > 2\lambda_3(2)$, there exists a constant $C_7 \equiv C_7(s, \kappa, \lambda) \in (0, \infty)$ such that for any $n \geq 2$,*

$$E[(1 + [\sqrt{b}|\hat{\Delta}_n|]^{s_0})\mathbb{1}(\sqrt{b}|\hat{\Delta}_n| > [(s-2)\lambda \log n]^{1/2})] \leq C_7 b^{-(s-2)/2}(\log n)^{-2}.$$

3.2. *Expansions for the Studentized sample mean.* Let $\{X_i\}_{i\in\mathbb{Z}}$ be a *real-valued* (i.e., $d_0 = 1$) stationary process. As an application of the main results of Section 2, here we establish EEs for a version of the Studentized sample mean, given by

$$(3.4) \qquad T_n \equiv n^{1/2}(\bar{X}_n - EX_1)/\tilde{\sigma}_n,$$

where $\tilde{\sigma}_n^2 = \max\{n^{-1}, b^{-1}\sum_{i=1}^{b}[U_{1(i-1)\ell+1} - (b^{-1}\sum_{i=1}^{b}U_{1(i-1)\ell+1})]^2\}$ is the nonoverlapping block bootstrap (NBB) (cf. Carlstein [4]) estimator of $n\operatorname{Var}(\bar{X}_n)$ based on blocks of length $\ell$, truncated from below at $n^{-1}$ and where $b = n/\ell$ (assumed to be an integer, for simplicity). We also assume that

$$(3.5) \qquad \ell \sim \beta_0 n^{1/3} \qquad \text{for some } \beta_0 \in (0, \infty),$$

where, for $\{u_n\}_{n\geq 1}, \{v_n\}_{n\geq 1} \in (0, \infty)$, we write $u_n \sim v_n$ if $\lim_{n\to\infty} u_n/v_n = 1$. Note that (3.5) covers the optimal block size for estimating $\operatorname{Var}(\bar{X}_n)$ by the NBB (cf. Hall, Horowitz and Jing [12], Lahiri [20]). Since the block variables in this problem are smooth functions of the $X_i$'s, we may impose regularity conditions on the $X_i$'s directly, without any reference to the block variables.

CONDITIONS.

S.1. $\{X_i\}_{i\in\mathbb{Z}}$ is stationary, $\rho \equiv E|X_1|^{2(s+1)}\{\log(1+|X_1|)\}^{\gamma(s)} \in (0, \infty)$ for some $s \geq 3$, $s \in \mathbb{N}$ and $\gamma(s) > s^2$, $EX_1 = 0$ and $\sigma_\infty^2 = \sum_{i\in\mathbb{Z}} EX_0X_i \in (0, \infty)$.
S.2. (i) The $\sigma$-fields $\mathcal{D}_j$ are generated by a sequence of *independent* $d_2$-dimensional ($d_2 \in \mathbb{N}$) random vectors $Z_j$, that is, $\mathcal{D}_j = \sigma\langle Z_j\rangle, j \in \mathbb{Z}$.

   (ii) There exists $\kappa \in (0,1)$ such that for all $m > \kappa^{-1}$ and $j \geq 1$, there exists a $\mathcal{D}_{j-m}^{j+m}$-measurable $X_{j,m}^\dagger$ such that $E\|X_j - X_{j,m}^\dagger\|^3 \leq \kappa^{-1}\exp(-\kappa m)$.
   (iii) There exists $\kappa \in (0,1)$ such that for all $m, j_0 \in \mathbb{N}$ with $\kappa^{-1} < m < j_0$,

$$\sup\{E|E\{\exp(\iota t[X_{j_0-m} + \cdots + X_{j_0+m}])|\mathcal{D}_j : j \neq j_0\}| : |t| \geq \kappa\} \leq 1 - \kappa.$$

Condition S.2(iii) is the conditional Cramér condition introduced by [8] for deriving valid asymptotic expansions for $S_{1n}$. Validity of a general order EE for the Studentized mean $T_n$ under S.2(iii) is somewhat *surprising* as the validity of a similar expansion for the Studentized sample mean for i.i.d. variables requires a stronger version of the Cramér condition that involves both the linear and the quadratic parts. This phenomenon can be explained by noting that in the dependent case, the block variables $U_{1i}$ tend to behave like their limit, that is, like $N(0, \sigma_\infty^2)$ variates as the block size $\ell \to \infty$ and the (stronger) Cramér condition for the *bivariate* random vector $(Z_\infty, Z_\infty^2)'$ holds when $Z_\infty$ is a $N(0, \sigma_\infty^2)$ variate. We now describe two important classes of weakly dependent processes for which condition S.2 holds. See [8] and [9] for more examples.

EXAMPLE 3.1 (*Linear processes*). Let $X_j = \sum_{i \in \mathbb{Z}} a_i \varepsilon_{j-i}, j \in \mathbb{Z}$, where $\{\varepsilon_i\}_{i \in \mathbb{Z}}$ are i.i.d. random variables with $E\varepsilon_1 = 0$ and $\sigma^2 = E\varepsilon_1^2 \in (0, \infty)$, and where $\{a_i\}_{i \in \mathbb{Z}} \in \mathbb{R}$ is such that $\sum_{i \in \mathbb{Z}} a_i \neq 0$ and for some $\kappa \in (0, \infty)$, $|a_i| = O(\exp(-\kappa|i|))$ as $|i| \to \infty$. Then S.2 holds with $\mathcal{D}_j = \sigma\langle \varepsilon_j \rangle, j \in \mathbb{Z}$, provided $E|\varepsilon_1|^3 < \infty$ and $\varepsilon_1$ satisfies the usual Cramér condition (cf. [8]) $\limsup_{|t| \to \infty} |E \exp(\iota t \varepsilon_1)| < 1$.

EXAMPLE 3.2 (*m-dependent processes*). Let $X_i = h_\varepsilon(\varepsilon_i, \ldots, \varepsilon_{i+m_0-1})$, $i \in \mathbb{Z}$, where $\{\varepsilon_i\}_{i \in \mathbb{Z}}$ are i.i.d. random variables with Lebesgue density $f_\varepsilon$, $m_0 \in \mathbb{N}$, and $h_\varepsilon \in C^1(\mathbb{R}^{m_0})$. Then, $\{X_i\}_{i \in \mathbb{Z}}$ is an $m_0$-*dependent* process. In this case, condition S.2 holds with $\mathcal{D}_j = \sigma\langle \varepsilon_j \rangle$, and $X_{j,m}^\dagger = X_j$ for $m \geq m_0$, $j \in \mathbb{Z}$, provided (cf. [8]) that there exist points $y_1, \ldots, y_{2m_0-1} \in \mathbb{R}$ and an open set $\mathcal{O} \subset \mathbb{R}$ such that $y_j \in \mathcal{O}$ for all $1 \leq j \leq 2m_0 - 1$, the density $f_\varepsilon$ is strictly positive on $\mathcal{O}$ and $0 \neq \sum_{j=1}^{m_0} \frac{\partial}{\partial x_j} h_\varepsilon(x_1, \ldots, x_{m_0})|_{(x_1,\ldots,x_{m_0})=(y_j,\ldots,y_{j+m_0-1})}$.

The following is the main result of this section. Write $\phi(x)$ and $\Phi(x)$ to denote the standard normal density and distribution function, respectively.

THEOREM 3.2. *Suppose that* (3.5) *and conditions* S.1 *and* S.2 *hold. Then there exist polynomials* $p_{rn}, 1 \leq r \leq s-2$, *with bounded coefficients such that*

$$\sup_{x \in \mathbb{R}} \left| P(T_n \leq x) - \left\{ \Phi(x) + \sum_{r=1}^{s-2} n^{-r/3} p_{rn}(x) \phi(x) \right\} \right| = O(n^{-(s-2)/3}(\log n)^{-2}).$$

The coefficients of $p_{rn}$ are $O(1)$ as $n \to \infty$ and typically contain smaller-order terms. For example, the third-order EE for $T_n$ (with $s = 5$) has the form

$$P(T_n \leq x) = \int_{-\infty}^x \phi(y)[1 + \{n^{-1/3}p_1(y) + n^{-1/2}p_2(y) + n^{-2/3}p_3(y)$$
$$+ n^{-5/6}p_4(y) + n^{-1}p_5(y)] \, dy$$
$$+ O(n^{-1}(\log n)^{-2})$$

uniformly in $x \in \mathbb{R}$, where the $p_j$'s are some polynomials whose coefficients are rational functions of the moments of the block variables $U_{1i}$. In particular, for $y \in \mathbb{R}$,

$$p_1(y) = (y^2 - 1)\left[\sum_{k=1}^\infty k \operatorname{Cov}(X_1, X_{k+1})\right][\sigma_\infty^{-2}][n^{1/3}/\ell],$$

$$p_2(y) = \left[-\frac{\ell^{1/2}EZ_nV_n}{2(EU_{11}^2)^{3/2}}\right]y$$
$$+\frac{1}{6}\left[(n^{1/2}EZ_n^3)(EU_{11}^2)^{-3/2} - \frac{3}{(EU_{11}^2)^{5/2}}(EZ_n^2)(\ell^{1/2}EZ_nV_n)\right]$$
$$\times (y^3 - 3y),$$

where $Z_n = n^{1/2}(\bar{X}_n - EX_1)$ and $V_n = b^{-1/2}\sum_{i=1}^{b}[U_{1i}^2 - EU_{1i}^2]$. As in the independent case (cf. Bhattacharya and Ghosh [2]), the terms in the EEs can be derived by "formally" expanding the cumulant generating function of the polynomial stochastic approximation to the Studentized statistic, and then by Fourier inversion of the resulting series. See Lahiri [21] for more details. Note that, as mentioned in Section 1, the expansion for the Studentized sample mean here is a combination of two series in powers of $n^{-1/3}$ and $n^{-1/2}$, the first corresponding to the (cumulants of the) Studentizing factor and the second coming from the expansion for $S_{1n}$.

3.3. *Second-order correctness of the blocks of blocks bootstrap.* For bootstrapping estimators of infinite-dimensional parameters of a weakly dependent process, Politis and Romano [24] formulated the BOBB method, (a version of) which is now described briefly. Given $\{X_1, \ldots, X_n\}$, let $Y_{in} = f_{in}(\chi_{i,\ell}), i = 1, \ldots, N$, be a set of block variables based on overlapping blocks $\chi_{i,\ell}$ of length $\ell$, where $N = n - \ell + 1$ and the $f_{in}$'s are Borel measurable functions [cf. (1.1)]. To approximate $\mathcal{L}(\check{T}_N)$, the distribution of $\check{T}_N \equiv t_N(Y_{1n}, \ldots, Y_{Nn}; \theta)$ using the BOBB, one defines the blocks of block variables $\{\mathcal{B}_i \equiv (Y_{in}, \ldots, Y_{(i+\ell_1-1)n}) : i = 1, \ldots, N - \ell_1 + 1\}$, resamples $\lceil N/\ell_1 \rceil$ many blocks from $\{\mathcal{B}_1, \ldots, \mathcal{B}_{N-\ell_1+1}\}$, and uses the first $N$ resampled block variables $\{Y_{in}^* := 1, \ldots, N\}$ (say) to define the BOBB version of $\check{T}_N$ as $\check{T}_N^* = t_N(Y_{1n}^*, \ldots, Y_{Nn}^*; \hat{\theta})$, where $\hat{\theta}_n$ is an estimator of $\theta$ based on $X_1, \ldots, X_n$. Then $\mathcal{L}(\check{T}_n^*)$ (conditional on $X_1, \ldots, X_n$) is the desired BOBB estimator of $\mathcal{L}(\check{T}_N)$.

Politis and Romano [24] established consistency of the BOBB method for statistics of the form $\tilde{T}_N \equiv b^{1/2}(\bar{Y}_N - E\bar{Y}_N)$, where $\bar{Y}_N = N^{-1}\sum_{i=1}^{N}Y_{in}$. We now show that the BOBB is indeed s.o.c. for a *Studentized* version of $\tilde{T}_N$, given by

(3.6) $$T_N \equiv \sqrt{b}(\bar{Y}_N - E\bar{Y}_N)/\hat{\sigma}_n,$$

where $\hat{\sigma}_n^2 = \max\{n^{-1}, [\hat{\gamma}_n^0(0) + 2\sum_{k=1}^{2\ell}(1 - N^{-1}k)\hat{\gamma}_n^0(k)]b/N\}$ and $\hat{\gamma}_n^0(k) = N^{-1}\sum_{i=1}^{N-k}(Y_{in} - \bar{Y}_N)(Y_{(i+k)n} - \bar{Y}_N)$, $k = 0, \ldots, N-1$. To define the BOBB version of $T_N$, let $\bar{Y}_{in,1}^* = \ell_1^{-1}\sum_{j=((i-1)\ell_1+1}^{i\ell_1}Y_{in}^*$, $i = 1, \ldots, b_1$, where $b_1 = N/\ell_1$, which for simplicity of exposition is assumed to be an integer. Then the BOBB version of $\bar{Y}_N$ is given by $\bar{Y}_N^* \equiv b_1^{-1}\sum_{i=1}^{b_1}\bar{Y}_{in,1}^*$. And for $\hat{\sigma}_n^2$, we use a truncated version of the sample variance of the $\bar{Y}_{in,1}^*$'s, namely, $[\sigma_n^*]^2 \equiv \max\{n^{-1}, b_1^{-1}\sum_{i=1}^{b_1}[\bar{Y}_{in,1}^* - \bar{Y}_N^*]^2\}$, to define the BOBB version of $T_N$ as (cf. Götze and Künsch [10])

$$(3.7) \qquad T_N^* = \sqrt{b_1}(\bar{Y}_N^* - E_*\bar{Y}_N^*)/\sigma_n^*.$$

Since the rate of normal approximation to $T_N$ is $O(b^{-1/2})$, here s.o.c. refers to an error of approximation that is of smaller order than $O(b^{-1/2})$. The next result gives sufficient conditions for the s.o.c. of $T_N$.

THEOREM 3.3. *Suppose that $\{X_i\}_{i\in\mathbb{Z}}$ is a sequence of stationary $\mathbb{R}^{d_0}$-valued random vectors, $Y_{in} = f_n(\chi_{i,\ell})$ for some $f_n:\mathbb{R}^{d_0\ell} \to \mathbb{R}$. Also, suppose that:*

(i) *for some $\kappa \in (0,1)$ and $\delta_0 \in (0,1/5)$, $\ell < \kappa^{-1}n^{1/5}$ and $\ell_1 < \kappa^{-1}n^{\delta_0}\ell$ for all $n > \kappa^{-1}$, and $[\log n]/\ell + \ell/\ell_1 = o(1)$ as $n \to \infty$;*

(ii) *for some $\kappa \in (0,1)$ and for some integer $r_0 \geq 12[(1-5\delta_0)(4-5\delta_0)]^{-1}$ [where $\delta_0$ is as in (i) above], $\sup_{n\geq 1} E|Y_{1n}|^{4r_0+\kappa} < \infty$;*

(iii) $\sigma_\infty^2 \equiv \lim_{n\to\infty}\mathrm{Var}(b^{1/2}\bar{Y}_N)$ *exists, $\sigma_\infty^2 \in (0,\infty)$; and*

(iv) *conditions* C.4, C.5 *hold and condition* C.6 *holds with $d_n = 1$ and $A_n$ replaced by $A_n^0 \equiv\equiv \{(0', s)' : s \in \mathbb{R}, |s| \geq \kappa\}$. Then*

$$(3.8) \qquad \sup_{x\in\mathbb{R}} |P(T_N \leq x) - P_*(T_N^* \leq x)| = o_P(b^{-1/2}) \qquad as\ n \to \infty.$$

Theorem 3.3 shows that for s.o.c. of the BOBB, the BOBB block size $\ell_1$ must be of a larger order than $\ell$. In addition, we also need an upper bound on the growth rate of $\ell_1/\ell$. Since $T_N^*$ is defined in terms of $b_1$-many resampled blocks, its EE has an error of the order $o(b_1^{-1/2})$, which must be $o(b^{-1/2})$ for the BOBB to be s.o.c. A large $\ell_1$ can make $b_1$ too small compared to $b$ and s.o.c. cannot be attained. See [22] for more details. We also mention that the Cramér condition (iv) is only on the sum of the block variables $Y_{in}$, not on the joint distribution of $Y_{in}$ and $Y_{(i+k)n}Y_{in}$, $0 \leq k \leq \ell_1$. Under this setup, the proof heavily depends on the arguments developed in Lahiri [19], which established a similar s.o.c. result for the MBB method in the *finite*-dimensional parameter case.

For a concrete example where Theorem 3.3 can be readily applied, suppose that the $X_i$'s are real-valued and that $Y_{in} = (2\pi\ell)^{-1}|\sum_{j=i}^{i+\ell-1} X_j \exp(-\iota j\omega)|^2$, $\omega \in [-\pi,\pi]$. Then, $\bar{Y}_N$ gives an estimator of the spectral density $f$ of $\{X_i\}_{i\in\mathbb{Z}}$ and it is asymptotically equivalent to the estimator (1.3) of Example 1.1 with weights $\omega_{kn} \equiv 1$ (cf. Politis and Romano [24]). In this case, $\sigma_\infty^2$ exists and $\sigma_\infty^2 = (2/3)f^2(\omega)$ for $0 < |\omega| < \pi$ and $\sigma_\infty^2 = (4/3)f^2(\omega)$ for $\omega = 0, \pm\pi$. Thus, condition (iii) of Theorem 3.3 holds if $f(\omega) \in (0,\infty)$. If, in addition, we suppose that $\{X_i\}_{i\in\mathbb{Z}}$ is a linear process as in Example 3.1, then by the arguments of Janas [14], condition (iv) holds, provided $\mathcal{L}(\varepsilon_1)$ has an absolutely continuous component. Thus, for this class of linear processes, the BOBB approximation to the distribution of the Studentized spectral density estimator is s.o.c.

**4. Proofs.** For $x, y \in \mathbb{R}$, let $x \vee y = \max\{x, y\}$, $x_+ = x \vee 0$, and $\lfloor x \rfloor = k$ if $k \leq x < k+1$, $k \in \mathbb{Z}$. Let $m_3 = \lceil (\log n) \log \log(n+2) \rceil$. Let $\alpha! = \prod_{j=1}^d \alpha_j!$ and $t^\alpha = \prod_{j=1}^d t_j^{\alpha_j}$ for $\alpha = (\alpha_1, \ldots, \alpha_d)' \in (\mathbb{Z}_+)^d$, $t = (t_1, \ldots, t_d)' \in \mathbb{R}^d$. Let $\mathcal{F}_a^c = \sigma \langle X_i : i \in \mathbb{Z}, a \leq i \leq c \rangle$ for $-\infty \leq a \leq c \leq \infty$. Let $C_k, C$ and $C(\cdot)$ denote generic constants, depending on their arguments (if any), but not on $n$ and $\ell$; dependence on $d_0, d_1, d, \kappa, \rho$ will be often suppressed to ease notation. For $c > 0$, define the truncation function $g(\cdot; c) : \mathbb{R}^k \to \mathbb{R}^k$ by $g(x; c) = \frac{cx}{\|x\|} \cdot \psi(\frac{\|x\|}{c}) \mathbb{1}(x \neq 0)$, where $\psi \in C^\infty(\mathbb{R})$ is nondecreasing, satisfying $\psi(u) = u_+$ for $u \leq 1$ and $\psi(u) = 2$ for $u \geq 2$. Let

$$_1W_{jn} = \sqrt{\ell} \bar{X}_{jn}, \; _2W_{jn} = \bar{Y}_{jn} \quad \text{and} \quad W_{jn} = (_1W_{jn}'; _2W_{jn}')', \qquad 1 \leq j \leq b.$$

With $c_n = b^{1/2} (\log n)^{-2}$, define the truncated and the centered versions of $_k W_{jn}$ by $_k \check{W}_{jn} = g([_k \check{W}_{jn}]; c_n)$, $_k \tilde{W}_{jn} = [_k \check{W}_{jn}] - E[_k \check{W}_{jn}]$ and set $_k \check{S}_n = \tilde{b}^{-1/2} \sum_{j=1}^b [_k \check{W}_{jn}]$ and $_k \tilde{S}_n = \tilde{b}^{-1/2} \sum_{j=1}^b [_k \tilde{W}_{jn}]$, $k = 1, 2$. Let $\check{W}_{jn} = ([_1 \check{W}_{jn}]', [_2 \check{W}_{jn}]')'$ and similarly, define $\tilde{W}_{jn}, \check{S}_n$ and $\tilde{S}_n$. For any random variable (r.v.) $Z$ on $(\Omega, \mathcal{F}, P)$, let

$$E_t Z = E[Z \exp(\iota t' \tilde{S}_n)] / H_n(t), \qquad t \in \mathbb{R}^d,$$

where $H_n(t) = E \exp(\iota t' \tilde{S}_n)$. Define the semiinvariants of the r.v.s $V_1, \ldots, V_p$ by

$$(4.1) \quad \begin{aligned} &\mathcal{K}_t(V_1, \ldots, V_p) \\ &= \frac{\partial}{\partial \varepsilon_1} \cdots \frac{\partial}{\partial \varepsilon_p} \Big|_{\varepsilon_1 = \cdots = \varepsilon_p = 0} \log E \exp(\iota t' \tilde{S}_n + \varepsilon_1 V_1 + \cdots + \varepsilon_p V_p), \end{aligned}$$

$t \in \mathbb{R}^d$. Write $\mathcal{K}_t(V_1^p, V_2^q) = \mathcal{K}_t(V_1, \ldots, V_1, V_2, \ldots, V_2)$ where, on the right-hand side, $V_1$ appears $p$ times and $V_2$ appears $q$ times. Then, by Taylor's expansion, we get

$$(4.2) \qquad \log H_n(t) = \sum_{r=2}^s \mathcal{K}_0(t' \tilde{S}_n^r) + R_n^*(t),$$

where $R_n^*(t) = [\int_0^1 (1 - \eta)^s \mathcal{K}_{\eta t}(t' \tilde{S}_n^{(s+1)}) d\eta] / s!$.

LEMMA 4.1. *Let conditions* C.1–C.4 *hold. Then there exists* $C = C(\rho, s, \kappa)$ *such that for any* $a_1, \ldots, a_r \in \mathbb{R}^d$ *with* $\|a_i\| = 1$, $2 \leq r \leq s$, *and for all* $n \geq \kappa^{-1}$,

(i) $|\mathcal{K}_0(a_1'S_n,\ldots,a_r'S_n) - \mathcal{K}_0(a_1'\tilde{S}_n,\ldots,a_r'\tilde{S}_n)| \leq Cb^{-(s-2)/2} \times (\log n)^{2s-2-[\alpha(s)/s]}$,

(ii) $|\mathcal{K}_0(a_1'\tilde{S}_n,\ldots,a_r'\tilde{S}_n)| \leq Cb^{-(r-2)/2}$.

PROOF. The left-hand side of the inequality in part (i) is bounded above by

$$(4.3) \qquad \tilde{b}^{-r/2} \sum_{i=0}^{b-1} \sum^{(i,r)} |\mathcal{K}_0(V_{j_1},\ldots,V_{j_r}) - \mathcal{K}_0(\tilde{V}_{j_1},\ldots,\tilde{V}_{j_r})|,$$

where the summation $\sum^{(i,r)}$ extends over all $1 \leq j_1 \leq \cdots \leq j_r \leq b$ with maximal gap $i$, and where $V_{j_p} = a_p' W_{j_p n}$ and $\tilde{V}_{j_p} = a_p' \tilde{W}_{j_p n}, 1 \leq p \leq r$. Since $W_{in} \in \mathcal{F}_{(i-1)\ell+1}^{(i+1)\ell}$, $W_{i+j,n} \in \mathcal{F}_{(i+j-1)\ell+1}^{\infty}$ and $[(i+j-1)\ell - (i+1)\ell] = (j-2)\ell$, by C.3 and C.4 above and by (3.13) and (3.14) of [8], we have (see [22])

$$(4.4) \qquad \sum_{i \geq C(\kappa)} \sum^{(i,r)} |\mathcal{K}_0(\tilde{V}_{j_1},\ldots,\tilde{V}_{j_r})| \leq C(\kappa,r)|J_2|b^r c_n^r \ell \exp(-3s\log n).$$

Truncating the $W_{jn}$'s at $e_i = \exp(\kappa i)$ and using C.3 and C.4, we get (cf. [22])

$$\sum_{i \geq i_0} \sum^{(i,r)} |\mathcal{K}_0(V_{j_1},\ldots,V_{j_r})|$$

$$\leq \sum_{i \geq i_0} \sum^{(i,r)} |\mathcal{K}_0(a_1' g(W_{j_1 n}; e_i),\ldots,a_r' g(W_{j_r n}; e_i))|$$

$$(4.5) \qquad + \sum_{i=i_0}^{\infty} \sum^{(i,r)} C(r,\rho) \max\left\{ E\left(\prod_{p \in I} \|W_{j_p n}\| \mathbb{1}(\|W_{j_k n}\| > e_i)\right) : k \in I, \right.$$

$$\left. 1 \leq |I| \leq r \right\}$$

$$\leq C(\rho,\kappa,r,s) i_0^{r-\alpha(s)/s} b$$

for all $2 \leq r \leq s$, where $i_0 = \lceil \log n \rceil$. Also, by similar arguments, for $2 \leq r \leq s$,

$$(4.6) \qquad \sum_{i < i_0} \sum^{(i,r)} |\mathcal{K}_0(V_{j_1},\ldots,V_{j_r}) - \mathcal{K}_0(\tilde{V}_{j_1},\ldots,\tilde{V}_{j_r})|$$

$$\leq Cb(\log n)^r / \{c_n^{(s-r)}[\log n]^{\alpha(s)/s}\}.$$

Hence, part (i) follows from (4.3)–(4.6). Proof of part (ii) is similar; see [22]. □

For the next lemma, let $\tilde{W}_n(I) \equiv \prod_{j \in I} \prod_{p=1}^{r_j} a'_{jp} \tilde{W}_{jn}$ and $S_I^{(r)} \equiv S_I^{(r)}(m) = \iota \tilde{b}^{-1/2} t' \sum^{*r} \tilde{W}_{jn}$, $I \subset \{1,\ldots,b\}$, $r \geq 0$, $m \geq 3$, where $a_{jp} \in \mathbb{R}^d$ with $\|a_{jp}\| = 1$, $r_j \in \mathbb{N}$, and where $\sum^{*r}$ extends over all $1 \leq j \leq b$ with $|j - k| > mr$ for all $k \in I$.

LEMMA 4.2. *Let conditions* C.1–C.4 *hold. Then, for any* $\eta \in (0, \frac{1}{4})$, *there exists* $C(\eta, \rho) > 0$ *such that for every* $3 \leq m \leq \eta b/|I|$ *and* $1 \leq K \leq \eta b/(m|I|)$,

$$（4.7） \quad |H_n(t)||E_t \tilde{W}_n(I)| \leq |E\tilde{W}_n(I)|[\max\{|E \exp(S_I^{(r)})| : 1 \leq r \leq K\} + \eta^K] + C(|I|) b c_n^\gamma K 2^K \ell \exp(-C(\kappa) m\ell)$$

*for all* $\|t\| < C(\eta, \rho)[b/m]^{1/2}$, *where* $\gamma = \sum_{j \in I} r_j$.

PROOF. Let $\Delta_{1,r} = [\exp(S_I^{(r-1)} - S_I^{(r)}) - 1]$ and $\tilde{W}_{1,I} = \tilde{W}_n(I) \exp(\iota \tilde{b}^{-1/2} \times \sum_{j \in I} t' \tilde{W}_{jn})$, $r \geq 1$. Using Tikhomirov's [29] iterative method, one can show that the left-hand side of (4.7) $\leq \sum_{r=1}^{K} |E\tilde{W}_{1,I}(\prod_{j=1}^{r-1} \Delta_{1,j}) \exp(S_I^{(r)})| + |E\tilde{W}_{1,I}(\prod_{j=1}^{K} \Delta_{1,j}) \exp(S_I^{(K)})|$. Now, approximating $\Delta_{1,j}$'s and $S_I^{(r)}$'s using $X_{j,m_0}^\dagger$'s and $Y_{jn,m_0}^\dagger$'s (with $m_0 = \lfloor m\ell/12 \rfloor$), and using conditions C.3, C.4 and the bound "$\max\{E(\sum_{j=j_0+1}^{j_0+m} \tilde{W}_{jn})^2 : 1 \leq j_0 \leq b - m\} \leq C(\rho)m$, $m \geq 3$", one can complete the proof. See [22] for details. □

LEMMA 4.3. *Let* C.1–C.4 *hold and let* $I_1, I_2 \subset \{1, \ldots, b\}$ *with* $\min\{I_2\} - \max\{I_1\} \geq m_1$ *for some* $3 \leq m_1 \leq b - |I_1| - |I_2|$. *Then, given* $\eta \in (0, \frac{1}{4})$ *and an integer* $m \in [3, m_1/3]$, *there exists* $C_1 = C(\rho, \kappa, \eta) \in (0, \infty)$ *such that*

$$（4.8） \quad |E_t \tilde{W}_n(I_1) \tilde{W}_n(I_2) - E_t \tilde{W}_n(I_1) E_t \tilde{W}_n(I_2)|$$
$$\leq C(r, k, \eta) b c_n^\gamma [K 2^K \ell \exp(-C(\kappa) m\ell) + \eta^K] |H_n(t)|^{-2}$$

*for all* $\|t\| < C_1[b/m]^{1/2}$, $1 \leq K \leq m_1/m$, *and* $n \geq 1$, *where* $\gamma = \sum_{p=1}^{2} \sum_{j \in I_p} r_j$.

PROOF. One can prove (4.8) by using Tikhomirov's [29] method and arguments in the proofs of Lemma 3.1 of [17] and Lemma 4.2 above. See [22] for details. □

LEMMA 4.4. *Let conditions* C.1–C.4 *hold. Given* $\eta \in (0, \frac{1}{4})$ *and* $q \in \mathbb{Z}_+$,

$$（4.9） \quad \left|\frac{\partial^q}{\partial u^q} R_n^*(t + x_0 u)\right|_{u=0} \leq \frac{C(s, q, \rho, \kappa, \eta)(1 + \beta_n(t)^p)(1 + \|t\|^p)}{b^{(s-2)/2}(\log n)^{[\alpha(s) - 2s - 1]}}$$

*for all* $\|x_0\| \leq 1$ *and for all* $t \in A_n \equiv \{x \in \mathbb{R}^d : \|x\| < C(\eta, \rho, \kappa) b^{(1-\eta)/2}, \beta_n(x) < \infty\}$, *where* $\beta_n(t) = |H_n(t)|^{-1}[\sup\{|E \exp(S_I^{(r)})| : 1 \leq r \leq C(p_0) b^{1-\eta}, |I| \leq p_0\} + \exp(-C(\eta, \kappa, p_0) m_3 \ell)]$, $p_0 = s + 1 + q$, *and* $S_I^{(r)}$ *is as in Lemma* 4.2.

PROOF. It is enough to consider $b^{-r/2}\sum_{i=0}^{b-1}\sum^{(i,r)}|\mathcal{K}_{ut}(V_{j_1},\ldots,V_{j_r})|$ for $s+1 \leq r \leq p_0$ and $0 \leq u \leq 1$ [cf. (4.2)], where $V_j = \iota t' \tilde{W}_{jn}$ and $\sum^{(i,r)}$ is as in (4.3). By Lemma 4.1, for any $s+1 \leq r \leq p_0$, $E\|\tilde{W}_{jn}\|^r \leq c_n^{r-s-1} E\|\tilde{W}_{jn}\|^{s+1} \leq C(\rho,s) c_n^{r-s} (\log n)^{-\alpha(s)}$, uniformly in $1 \leq j \leq b$. Next, define $\Gamma_1 = \{i \in \mathbb{Z} : 1 \leq i \leq a_n\}$, $\Gamma_2 = \{i \in \mathbb{Z} : a_n < i \leq b^\eta\}$ and $\Gamma_3 = \{i \in \mathbb{Z} : b^\eta < i \leq b-1\}$, where $a_n = m_3^2$. Now using Lemma 4.2 (with $m = \lceil b^\eta \rceil$ and $K = m$) for $i \in \Gamma_1$, Lemma 4.3 (with $m = K = \lfloor \sqrt{i} \rfloor + 1$) for $i \in \Gamma_2$ and again Lemma 4.3 (with $m = K = \lfloor b^\eta \rfloor + 1$) for $i \in \Gamma_3$, one can complete the proof of the lemma. See [22] for more details. □

LEMMA 4.5. *Suppose that conditions* C.1–C.5 *hold. Then, for any* $I \subset \{1,\ldots,b\}$ *with* $|I| \equiv r \leq C$, *and for any* $t \in \mathbb{R}^d$, $3 \leq m \leq b/C$,

$$|H_n(t) E_t \tilde{W}_n(I)| \leq C c_n^\gamma [\beta_{1n}(t)]^K + C(\kappa,r) K c_n^\gamma [1 + \|t\|\ell m] \exp(-C(\kappa)m\ell)$$

*for some* $K \geq C[b/m]$, *where* $\beta_{1n}(t) = \max\{E|E(\exp(\iota \tilde{b}^{-1/2} t' \sum_{|j-j_0| \leq m} \tilde{W}_{jn})| \tilde{\mathcal{D}}_{j_0})| : m < j_0 < b-m\}$, $\gamma = \sum_{j \in I} r_j$ *and* $\tilde{\mathcal{D}}_{j_0}$ *is as in condition* C.6.

PROOF. Let $I = \{j_1,\ldots,j_r\}$, $J_{0n} = \{1,\ldots,b\}$ and $J_{1n} = \{j \in J_{0n} : |j - j_k| \geq 2m+1 \text{ for all } 1 \leq k \leq r\}$. Define $j_1^0 = \inf J_{1n}$ and $j_{p+1}^0 = \inf\{j \in J_{1n} : j \geq j_p^0 + 7m\}$, $p = 1, 2, \ldots, K-1$, where $K$ is the first integer $p$ for which the infimum is over an empty set. Also, define the variables $A_p = \exp(\iota \tilde{b}^{-1/2} t' \sum_{|j-j_p^0| \leq m} \tilde{W}_{jn})$, $B_p = \exp(\iota \tilde{b}^{-1/2} t' \sum_{j_p^0+m+1 \leq j \leq j_{p+1}^0 - m - 1} \tilde{W}_{jn})$ and $R = \tilde{W}_n(I) \exp(\iota \tilde{b}^{-1/2} t' \times \sum_{j \in J_{2n}} \tilde{W}_{jn})$, where $J_{2n} = \{j \in J_{0n} : j < j_1^0 - m \text{ or } j \geq j_K^0 - m\}$. Then it follows that

$$H_n(t) E_t \tilde{W}_n(I) = E\left[\left(\prod_{p=2}^{K-1} A_p B_p\right) R\right].$$

Let $A_p^\dagger, B_p^\dagger$ and $R^\dagger$ be defined by replacing $X_j$'s and $Y_{jn}$'s by $X_{j,q}^\dagger, Y_{jn,q}^\dagger$ in $A_p, B_p$ and $R$ with $q = m\ell$ for $A_p^\dagger$ and with $q = m\ell/12$ for $B_p^\dagger$ and $R^\dagger$. Then, by C.2,

$$\left| ER \prod_{p=2}^{K-1} A_p B_p - ER^\dagger \prod_{p=2}^{K-1} A_p^\dagger B_p^\dagger \right|$$
(4.10)
$$\leq C(\kappa,|I|) K c_n^\gamma \|t\|\ell m \exp(-c(\kappa)m\ell).$$

Next, let $\tilde{\mathcal{D}}_p = \sigma\langle\{\mathcal{D}_j : j \in \mathbb{Z}, j \notin [c_p, d_p]\}\rangle$ and $\mathcal{D}_p^* = \sigma\langle\{\mathcal{D}_j : j \in [a_p - m\ell, c_p) \cup (d_p, b_p + m\ell]\}\rangle$, where $a_p = (j_p^0 - m)\ell + 1 - m\ell$, $b_p = (j_p^0 + m + 1)\ell + m\ell$, $c_p = (j_p^0 - \lfloor m/2 \rfloor)\ell + 1$ and $d_p = (j_p^0 + \lfloor m/2 \rfloor + 1)\ell$, $2 \leq p \leq K-1$. Then, by condition C.5, $\max_{p=2,\ldots,K-1}|E(A_p^\dagger|\tilde{\mathcal{D}}_p)-E(A_p^\dagger|\mathcal{D}_p^*)|\leq C(\kappa)\exp(-\kappa m\ell)$ and hence,

$$\left|ER^\dagger\prod_{p=2}^{K-1}A_p^\dagger B_p^\dagger - ER^\dagger\prod_{p=2}^{K-1}B_p^\dagger E(A_p^\dagger|\mathcal{D}_p^*)\right|$$

$$\leq \sum_{q=2}^{K-1}\left|ER^\dagger\left(\prod_{p=2}^{q-1}A_p^\dagger B_p^\dagger\right)B_q^\dagger[A_q^\dagger-E(A_q^\dagger|\tilde{\mathcal{D}}_q)]\prod_{p=q+1}^{K-1}B_p^\dagger E(A_p^\dagger|\mathcal{D}_p^*)\right|$$

(4.11)

$$+Cc_n^\gamma \sum_{q=2}^{K-1}E|E(A_q|\tilde{\mathcal{D}}_q)-E(A_q^\dagger|\mathcal{D}_q^*)|$$

$$\leq C(\kappa)c_n^\gamma K\exp(-\kappa m\ell),$$

since by construction, $R^\dagger(\prod_{p=2}^{q-1}A_p^\dagger B_p^\dagger), B_q^\dagger$ and $\prod_{p=q+1}^{K-1}B_p^\dagger E(A_p^\dagger|\mathcal{D}_p^*)$ are measurable w.r.t. $\tilde{\mathcal{D}}_q$ for every $2\leq q\leq K-1$, making the first term vanish. Now using (4.10), (4.11) and the fact that $\mathcal{D}_p^*$ and $\mathcal{D}_{p+1}^*$ are separated by a distance $\geq Cm\ell$ for all $p$, one can retrace the arguments in [8] to conclude that $|H_n(t)E_t\tilde{W}_n(I)|\leq Cc_n^\gamma\prod_{p=2}^{K-1}E|E(A_p|\tilde{\mathcal{D}}_p)|+C(\kappa,|I|)Kc_n^\gamma(1+\|t\|\ell m)\exp(-C(\kappa)m\ell)$. Lemma 4.5 follows from this. □

LEMMA 4.6. *Let* C.1–C.5 *hold and* $\eta\in(0,1)$. *Then, for any* $\alpha\in\mathbb{Z}_+^d$,

$$|D^\alpha H_n(t)|\leq C(\alpha,\eta,k)b^{|\alpha|}[\exp(-C(\kappa)\|t\|^2)+\exp(-C(\kappa)m_n\ell)]$$

*for all* $\|t\|\leq C(\kappa,\rho,s)b^{(1-\eta)/2}$, *where* $m_n=\lceil b^{\eta/2}\rceil$.

PROOF. For any $\mathbb{R}^d$-valued zero mean r.v. $Z$ on $(\Omega,\mathcal{F},P)$ with $E\|Z\|^3<\infty$ and sub-$\sigma$-field $\mathcal{C}\subset\mathcal{F}$, one can show (cf. [22]) that $|E\exp(\iota t'Z|\mathcal{C})|^2\leq 1-E((t'Z)^2|\mathcal{C})+2E(|t'Z|^3|\mathcal{C})$ for all $t\in\mathbb{R}^d$. Taking $Z=\tilde{b}^{-1/2}\sum_{j=j_0+1}^{j_0+m_0}\tilde{W}_{jn}$ and using C.2 and Hölder's inequality, one can show (cf. [22]) that

$$(4.12)\quad \left[E\left|E\left\{\exp\left(\iota t'\left[\tilde{b}^{-1/2}\sum_{j=j_0+1}^{j_0+m_0}\tilde{W}_{jn}\right]\right)\Big|\mathcal{C}\right\}\right|\right]^2\leq\exp\left(-\frac{\kappa}{4}b^{-1}m_0\|t\|^2\right)$$

for all $\|t\|\leq C(\kappa,\rho,s)b^{1/2}/m_0^2$, provided $n>\kappa^{-1}$. Now applying Lemma 4.5 with $m=\lceil b^{\eta/2}\rceil$ and using (4.12) to estimate $\beta_{1n}(t)$ of Lemma 4.4, one can complete the proof of Lemma 4.6. See [22] for more details. □

LEMMA 4.7. *Let conditions* C.1–C.4 *hold. Let* $f:\mathbb{R}^d\to\mathbb{R}$ *be a Borel measurable function with* $\sup\{\frac{|f(x)|}{(1+\|x\|^{s_0})}:x\in\mathbb{R}^d\}\equiv M_f<\infty$. *Then, for any* $\gamma>0$,

$$\left|Ef(S_n)-\int f\,d\Psi_{n,s}\right|$$

$$\leq C(\rho, s, \gamma, \Xi_\infty, M_f)$$

$$\times \left[ b^{-(s-2)/2}(\log n)^{2s-2-[\alpha(s)/s]} \right.$$

$$\left. + \sup_{|\alpha|\leq p_1} \int |D^\alpha[(H_n(t) - \hat{\Psi}_{s,n}(t))K_n(t)]|\,dt + \omega(\tilde{f}; b^{-\gamma}) \right],$$

where $p_1 = d + s_0 + 1$, $\tilde{f}(x) = \frac{f(x)}{1+\|x\|^{s_0}}$, $x \in \mathbb{R}^d$, $K_n(t) = K_0(\frac{t}{b^\gamma})$ and $K_0 \in C^{p_1}(\mathbb{R}^d)$ is a characteristic function that vanishes outside a compact set.

PROOF. Let $A_1 = \{\|S_n\| > \log n\}$, $A_2 = \{\|\check{S}_n\| > \log n\}$ and $A_3 = \{2\|\tilde{S}_n\| > \log n\}$, $n \geq 3$. Then, it can be shown (cf. [22]) that

$$I_1 \equiv |Ef(S_n) - Ef(\check{S}_n)|$$

(4.13)
$$\leq 12M_f[(\log n)^{s_0} P(S_n \neq \check{S}_n)$$
$$+ E\|\check{S}_n\|^{s_0} \mathbb{1}_{A_2} + |E\|S_n\|^{s_0} - E\|\check{S}_n\|^{s_0}|].$$

By Markov's inequality and Jensen's inequality, for all $n$ with $c_n > 1$, one gets

$$P(S_n \neq \check{S}_n) \leq \sum_{j=1}^{b} P(\|\sqrt{\ell}\bar{X}_{jn}\| > c_n) + \sum_{j=1}^{b} P(\|\bar{Y}_{jn}\| > c_n)$$

(4.14)
$$\leq \rho b[h_s(c_n)]^{-1} + [h_s(c_n)]^{-1} \sum_{j=1}^{b} E\left\{\ell^{-1} \sum_{i=(j-1)\ell+1}^{j\ell} h_s(\|Y_{in}\|)\right\}$$

$$\leq C(\rho, s) b^{-(s-2)/2}(\log n)^{-[\alpha(s)-2s]};$$

(4.15) $\quad \|E\check{S}_n\| \leq C(\rho, s) b^{-(s-2)/2}(\log n)^{2(s-1)-\alpha(s)}.$

Let $g_n(x) = \|x\|^{s_0}\mathbb{1}(2\|x\| > \log n)$ and $\tilde{g}_n(x) = g_n(x)(1 + \|x\|^{s_0})^{-1}$, $x \in \mathbb{R}^d$, $n \geq 2$. Then, for any $a > 0$, $\omega(\tilde{g}; b^{-a}) + \int g_n(x)\,d\Psi_{s,n}(x) \leq C(\rho, \Xi_\infty, a)n^{-a}$; see [22] for details. Hence, by Lemma 4.1 and by (4.13)–(4.15), for all $n > C$,

$$I_1 \leq C(\rho, \Xi_\infty, s) M_f\left[(\log n)^{2s-2-[\alpha(s)/s]} b^{-(s-2)/2} + \left|Eg_n(\tilde{S}_n) - \int g_n\,d\Psi_{s,n}\right|\right].$$

Lemma 4.7 now follows by two applications of the smoothing inequality of Sweeting [28] and Lemma 11.6 of Bhattacharya and Ranga Rao [3]. □

PROOF OF THEOREM 2.1. By Lemma 4.7, it is enough to show that for $0 \leq |\alpha| \leq p_1$,

(4.16)
$$\int |D^\alpha[(H_n(t) - \hat{\psi}_{s,n}(t))\hat{K}_0(b^{-a}t)]|\,dt$$
$$\leq C(d, s, \rho, \kappa) b^{-(s-2)/2}(\log n)^{-2}.$$

Suppose that $\hat{K}_0(t) = 0$ for all $\|t\| > C_0$. Then partition the set $\{\|t\| < C_0 b^a\}$ by $\{\|t\| \leq a_n\}$, $\{a_n < \|t\| < C(\kappa, \rho, s) b^{1-\eta}\}$ and $\{C(\kappa, \rho, s) b^{1-\eta} \leq \|t\| < C_0 b^a\}$, where $a_n = m_3^{1/2}$, $C(\kappa, \rho, s)$ is as in Lemma 4.6 and $\eta \in (0, \kappa/2)$. Then, using Lemma 4.4, one can show (cf. Lemma 3.33 of [8]) that for all $\alpha \in \mathbb{Z}_+^d$,

$$(4.17) \quad \int_{\|t\| < a_n} |[D^\alpha (H_n(t) - \hat{\Psi}_{s,n}(t)) \hat{K}_0(b^{-a} t)]| \, dt \leq \frac{C(s, |\alpha|, \rho, \kappa)}{b^{(s-2)/2} (\log n)^{[\alpha(s)-s-1]}}.$$

By condition C.2 and Lemma 4.6, for all $\alpha \in \mathbb{Z}_+^d$,

$$(4.18) \quad \begin{aligned} &\int_{\|t\| > a_n} |D^\alpha \hat{\Psi}_{s,n}(t)| \, dt + \int_{a_n < \|t\| < C(\kappa,\rho,s) b^{1-\eta}} |D^\alpha H_n(t)| \, dt \\ &\leq C(|\alpha|, \eta, \kappa) n^{-s}. \end{aligned}$$

Next, using condition C.6, one can show (cf. [22]) that for $n > C(\kappa)$,

$$\beta_{1n}(t) \leq \begin{cases} \exp(-C(\kappa)(d_n \tilde{b}^{1/2})^{-2} \|t\|^2), & \text{for all } \|t\| < \kappa \tilde{b}^{1/2} d_n, \\ \exp(-C(\kappa) \tilde{b}^{\kappa/2}), & \text{for all } b^{(1-\eta)/2} \leq \|t\| < \kappa \tilde{b}^{1/2} d_n. \end{cases}$$

Hence, by Lemma 4.5 (cf. the proof of Lemma 3.43 of [8]), one gets

$$(4.19) \quad \int_{C(\kappa, s\rho) b^{1-\eta} < \|t\| < C_0 b^a} |D^\alpha H_n(t)| \, dt \leq C(\rho, a, \kappa) n^{-s}.$$

Theorem 2.1 now follows from (4.16)–(4.19). □

PROOF OF THEOREM 2.2. Let $f_n(x) = f(x + E\check{S}_n)$ and $\eta \in (0, 1/s)$. Then a proof of Theorem 2.1 can be constructed by using the arguments in [7] and the expansion

$$f_n(x) = \sum_{0 \leq |\alpha| \leq s-2} (-\varepsilon)^{|\alpha|} t^\alpha D^\alpha f_n(x + \varepsilon t) \Big/ \left( \prod_{i=1}^d \alpha_i! \right)$$

$$+ (-\varepsilon)^{s-1} (s-1) \sum_{|\alpha|=s-1} \int_0^1 (1-u)^{s-2} D^\alpha f_n(x + u\varepsilon t) \, du \Big/ \left( \prod_{i=1}^d \alpha_i! \right),$$

with $\varepsilon = C(\{p(\alpha) : |\alpha| \leq s-1\}, \kappa, \rho) \cdot b^{(\eta-1)/2}$. See [22] for details. □

PROOF OF THEOREMS 2.3 AND 2.4. Use (4.17)–(4.18), Lemma 4.7, and the arguments in the proof of Theorem 2.11 of [8]. See [22] for more details. □

PROOF OF THEOREM 3.1. Let $Y_{jn} = (U_{1j}^\nu - EU_{1j}^\nu)$, $1 \leq j \leq N$ and $Y_{jn} = 0$ for $N+1 \leq j \leq n$. Then, by Theorem 2.4 and the stationarity of $\{Xi\}_{i \in \mathbb{Z}}$, it is enough to show that (i) $\limsup_{n \to \infty} \max_{1 \leq j \leq b} Eh_s(\bar{Y}_{jn}) < \infty$, (ii) condition C.3 holds and (iii) $\lim_{n \to \infty} \text{Var}(S_{2n})$ exists and lies in $(0, \infty)$.

By Lemma 3.3 of [17], $Eh_s(U_{11}^\nu) \leq C(|\nu|)(1 + E\|U_{11}\|^{2q}) = O(1)$, and hence, (i) holds. Next consider (ii). Let $Y_{jn,m}^\dagger = U_{1j,m}^\nu \mathbb{1}(\|U_{1j,m}\| \leq c_{1m}) - EU_{11}^\nu, 1 \leq j \leq N$ and $Y_{jn,m}^\dagger = 0$ for $N+1 \leq j \leq n$, where $U_{1j,m} = (X_{j,m}^\dagger + \cdots + X_{j+\ell-m}^\dagger)/\sqrt{\ell}$ and $c_{1m} = \exp(\kappa m/2|\nu|), m \geq 1$. Then, using the bound $E\|U_{11}\|^{2|\nu|} = O(1)$ (cf. Lemma 3.3 of [17]) and truncation arguments, one can show that

$$E|Y_{jn,m}^\dagger - Y_{jn}|$$
$$\leq E|U_{1j,m}^\nu - U_{1j}^\nu|\mathbb{1}(\|U_{1j,m}\| \leq c_{1m}) + E\|U_{1j}\|^r\mathbb{1}(\|U_{1j}\| > c_{1m})$$
$$\quad + E|U_{1j}^\nu||\mathbb{1}(\|U_{1j,m}\| \leq c_{1m}) - \mathbb{1}(\|U_{1j}\| \leq c_{1m})|$$
$$\leq C(\nu)[c_{1m}^{|\nu|-1}E\|U_{1j,m} - U_{1j}\| + E\|U_{1j}\|^{|\nu|}\mathbb{1}(\|U_{1j}\| > c_{1m}/2)]$$
$$\leq C(\nu,\kappa)\sqrt{\ell}\exp\left(-\frac{\kappa m}{2}\right)$$

for all $1 \leq j \leq N$, $n \geq 1, m \geq 1$; see [22] for details. Hence, (ii) holds.

Finally, consider (iii). Using (3.2), C.3 and C.4, one can show that $\sum_{j=\ell+3m}^{N-1}|\text{Cov}(U_{11}^\nu, U_{1(j+1)}^\nu)| = o(1)$, where $m = \lceil C(\kappa,s)\log n\rceil$. Next let $S(i,j) = X_i + \cdots + X_j, i \leq j$, $r = |\nu|$, and for $I \subset \{1,\ldots,r\}$, let $S_i(I) = \prod_{k \in I} t_k' S(1,i)$, where $t_1,\ldots,t_r$ are unit vectors such that $U_{1j}^\nu = \prod_{k=1}^r t_k' U_{1j}$. Then the stationarity of the $X_i$'s, (3.2), C.3 and C.4 yield (cf. [22])

$$EU_{11}^\nu U_{1(j+1)}^\nu = \ell^{-r}\sum_I \sum_J ES_{j-\ell^{1/4}}(I^c)ES_{j-\ell^{1/4}}(J^c)ES_{\ell-j}(I)S_{\ell-j}(J) + Q_\ell(j)$$

for $\ell^{1/2} \leq j \leq \ell - \ell^{1/2}$, where $\sum_I$ and $\sum_J$ run over $I, J \subset \{1,\ldots,r\}$ and where $\max\{|Q_\ell(j)|: \ell^{1/2} \leq j \leq \ell - \ell^{1/2}\} = o(1)$. Next using the CLT, uniform integrability of $S_{1n}^{2r}$ under (3.2), and the above expansion for $EU_{11}^\nu U_{1(j+1)}^\nu$, one can show that $\text{Var}(S_{2n}) = (n\ell)^{-1}\sum_{j=\ell^{1/2}}^{\ell-\ell^{1/2}}(N-j)\text{Cov}(U_{11}^\nu, U_{1(j+1)}^\nu) + o(1) = \lambda_3(\nu) + o(1)$ for some $\lambda_3(\nu) \in (0,\infty)$; see [22] for details. This completes the proof of (iii). □

PROOF OF COROLLARY 3.1. Without loss of generality (w.l.g.), we may set $\mu = 0$. Then, $b^{1/2}\hat{\Delta}_n = b^{1/2}[\hat{\mu}_n(2) - EU_{11}^2] - b^{1/2}[\hat{\mu}_n(2)]^2 + b^{1/2}[E(\hat{\mu}_n(2)]^2 \equiv T_{1n} + T_{2n} + b^{1/2}[E(\hat{\mu}_n(2)]^2$, say. Monotonicity of $h_c(x) = (1 + x^{s_0})\mathbb{1}(x > c), x \in [0,\infty)$ [for fixed $c \in (0,\infty)$] implies $Eh_c(|T_{1n}+T_{2n}|) \leq Eh_c(2\max\{|T_{1n}|, |T_{2n}|\}) \leq Eh_c(2|T_{1n}|) + Eh_c(2|T_{2n}|)$. Let $\varepsilon_0 = (2\lambda_3(2) - \lambda)/4$. Then using the bound $[E(\hat{\mu}_n(1)]^2 = O(n^{-1})$ (cf. Lahiri [20]) and applying Theorem 3.1 to $Eh_c(|T_{kn}|)$ with $c = [2\{\lambda_3(2) + \varepsilon_0\}(s-2)\log n]^{1/2}$ for $k = 1,2$, one can prove the result. □

PROOF OF THEOREM 3.2. Let $\sigma_\ell^2 = E(\sqrt{\ell}\bar{X}_2)^2$, $Y_{in}^{[2]} = (U_{1j}^{[2]})^2 - E(U_{11}^{([2])})^2$, $1 \leq j \leq b$, $S_{2n}^{[2]} = b^{-1/2}\sum_{j=1}^b Y_{jn}^{[2]}$ and $S_n^{[2]} = (n^{1/2}\bar{X}_n, S_{2n}^{[2]})'$, where $U_{1j}^{[2]} =$ $(X_{(j-1)\ell+1} + \cdots + X_{j\ell})/\sqrt{\ell}$. Then $T_n = \sqrt{n}\bar{X}_n/[\sigma_\ell^2 + b^{-1/2}S_{2n}^{[2]} - b^{-1}(\sqrt{n}\bar{X}_n)^2]$ is a smooth function of $S_n^{[2]}$ and hence, the EE for $T_n$ can be derived from that of $S_n^{[2]}$ (cf. Bhattacharya and Ghosh [2]). To derive an $(s-2)$th order EE for $S_n^{[2]}$, note that by (3.5) and the independence of $\{\mathcal{D}_j : j \in \mathbb{Z}\}$, it is enough to verify (cf. [22])

$$(4.20) \quad \lim_{n\to\infty} \text{Cov}(S_n^{[2]}) = \Xi_\infty^{[2]} \quad \text{exists and is nonsingular;}$$

$$(4.21) \quad \max_{j_0 \in J_n} \sup_{t \in A_n} E\left|E\left\{\exp\left(\iota t' \sum_{j=j_0-m_n}^{j_0+m_n} W_{jn}^{[2]}\right)\Big|\tilde{D}_{j_0}\right\}\right| \leq 1 - \kappa,$$

for some $\kappa \in (0,1)$, where $J_n, \tilde{D}_{j_0}$ are as in C.6, $m_n = \lceil C(\kappa,s)\log n \rceil$, $A_n = \{t \in \mathbb{R}^2 : \kappa \leq \|t\| \leq b^{(s-2)}\}$ and $W_{jn}^{[2]}$ is defined by replacing the $Y_{in}$'s in $W_{jn}$ with $Y_{in}^{[2]}$'s. Using S.1 and S.2(i)–(ii), one can show (cf. [22]) that (4.20) holds with $\Xi_\infty^{[2]} = \text{Diag}(\sigma_\infty^2, 2\sigma_\infty^4)$.

Next consider (4.21). Set $m_n = m$ for notational simplicity. Let $U_{2j,m}^\dagger = (X_{(j-1)\ell+1,m}^\dagger + \cdots + X_{j\ell,m}^\dagger)/\sqrt{\ell}$ and $c_{1m} = \exp(\kappa m/4)$. Then it can be shown that uniformly over $t \in A_n$ and $j_0 \in J_n$,

$$(4.22) \quad \sup_{t \in A_n} E\left|E\left\{\exp\left(\iota t' \sum_{j=j_0-m}^{j_0+m} W_{jn}^{[2]}\right)\Big|\tilde{\mathcal{D}}_{j_0}\right\}\right|$$
$$\leq E|E\{\exp(\iota(U_{2j_0,m}^\dagger, [U_{2j_0,m}^\dagger]^2)t)|\mathcal{D}_{j_0}^*\}| + 2P(\|U_{2j_0,m}^\dagger\| > c_{1m}),$$

where $\mathcal{D}_{j_0}^* \equiv \sigma\langle\{\mathcal{D}_j : j \in \mathbb{Z}, j \notin [(j_0-1)\ell+m+1, j_0\ell-m]\}\rangle \subset \tilde{\mathcal{D}}_{j_0}$. Next using the independence of the $\mathcal{D}_j$'s, one gets $E|E\{\exp(\iota(U_{2j,m}^\dagger,[U_{2j,m}^\dagger]^2)t)|\mathcal{D}_j^*\}| = E|\xi_{j,m}(t; Z(I_{2j}))|$ for some function $\xi_{j,m}(t;z)$, where for $I \subset \mathbb{Z}$, $Z(I) = \{Z(i) : i \in I\}$ and where $I_{2j} = \{i \in \mathbb{Z} : (j-1)\ell - m + 1 \leq i \leq (j-1)\ell + m\} \cup \{i \in \mathbb{Z} : j\ell - m + 1 \leq i \leq j\ell + m\}$. Let $Z_\infty \sim N(0, \sigma_\infty^2)$. Then, for any $\varepsilon > 0$, there exists $C(\varepsilon) \in (0,\infty)$ such that

$$(4.23) \quad \sup\{|E\exp(\iota(Z_\infty, Z_\infty^2)t)| : \|t\| > \varepsilon\} < e^{-C(\varepsilon)}.$$

Now using (4.23) and the EE results of [8], one can show that for any $\varepsilon > 0$, $\sup\{E|\xi_{j,m}(t; Z(I_{2j}))| : \varepsilon \leq \|t\| \leq b^{(s-2)}, j \in J_n\} \leq \exp(-C(\varepsilon)) + o(1)$; see [22] for details. Condition (4.21) now follows from (4.22), (4.23) and the above bound. □

PROOF OF THEOREM 3.3. W.l.g., let $EY_{1n} = 0$ for all $n \in \mathbb{N}$. For $k = 0, \ldots, N-1$, let $\tilde{g}_n^0(k) = N^{-1}\sum_{i=1}^{N-k} Y_{in}Y_{(i+k)n}$. Define $\tilde{\sigma}_n^2 = [\tilde{g}_n^0(0) + $ $2\sum_{k=1}^{2\ell}(1-N^{-1}k)\tilde{g}_n^0(k)]b/N$ and $\sigma_n^2 = E\tilde{\sigma}_n^2$, where we let $Z_N = b^{1/2}(\bar{Y}_N - E\bar{Y}_N)/\sigma_n$. Then by Taylor's expansion,

(4.24) $\quad T_N = Z_N - Z_N[\tilde{\sigma}_n^2 - \sigma_n^2]/[2\sigma_n^2] + R_{1N} \equiv T_{1N} + R_{1N},\qquad$ say,

where, on the set $\{|\hat{\sigma}_n^2 - \sigma_n^2| < \sigma_n^2/2\}$, $|R_{1N}| \leq \frac{3}{\sqrt{2}}\sigma_n^{-5}[\hat{\sigma}_n^2 - \sigma_n^2]^2 + |Z_n(\hat{\sigma}_n^2 - \tilde{\sigma}_n^2)|/[2\sigma_n^2]$. Next let $\eta_{1in} = Y_{in}[Y_{in} + 2\sum_{k=1}^{[2\ell]\wedge[N-i]}(1-N^{-1}k)Y_{(i+k)n}]$, $\eta_{in} = \eta_{1in} - E\eta_{1in}$, $i = 1,\ldots,N$ and $V_i(4\ell) \equiv \tilde{S}([(i-1)4\ell+1],[i4\ell \wedge N])$, where $\tilde{S}(p,q) = \sum_{i=p}^q \eta_{in}$ for any $1 \leq p \leq q \leq N$. Then, using the weak dependence of alternate $V_i(4\ell)$'s and Markov's inequality, with $\varepsilon_n = b^{-1/4}(\log n)^{-2}$, one gets

$$P(|\tilde{\sigma}_n^2 - \sigma_n^2| > 2\varepsilon_n)$$

(4.25)
$$\leq P\left(\left|\sum_{1\leq 2i\leq N/[4\ell]}V_{2i}(4\ell)\right| > \varepsilon_n N^2/b\right)$$
$$+ P\left(\left|\sum_{1\leq 2i-1\leq N/[4\ell]}V_{2i-1}(4\ell)\right| > \varepsilon_n N^2/b\right)$$
$$\leq C[\varepsilon_n N^2/b]^{-4}[\ell^4(N/[4\ell])^2] = O(b^{-1}(\log n)^8).$$

Also, by moderate deviation bounds for $\bar{Y}_N$ and Markov's inequality, $P(|\hat{\sigma}_n^2 - \tilde{\sigma}_n^2| > Cb^{-1/2}(\log n)^{-2}) = O(b^{-1/2}(\log n)^{-2})$ and hence, by (4.24)–(4.25), $T_N$ and $T_{1N}$ have identical EEs up to order $o(b^{1/2})$. With $a \equiv a_n = \lfloor b^{1/2}/m_3 \rfloor$, write $Z_N = Z_N^{(1)} + Z_N^{(2)}$ and $\tilde{\sigma}_n^2 - \sigma_n^2 = s_n^{(1)} + s_n^{(2)}$, where $Z_N^{(1)}$ denotes the sum over $Y_{in}b^{-1/2}$ for $i = 1,\ldots,2a\ell$ and $s_n^{(1)}$ over $N^{-2}b\eta_{in}, i = 1,\ldots,2a\ell$. Then it can be shown that

(4.26) $\quad T_{1N} = Z_N - Z_N^{(2)}s_n^{(2)}/[2\sigma_n^2] + R_{2N} \equiv T_{2N} + R_{2N},\qquad$ say,

where $P(|R_{2N}| > b^{-1/2}m_3^{-1/4}) = O(am_3^{1/2}/b) = o(b^{-/2})$. Thus, $T_N$ and $T_{2N}$ have identical EEs up to order $o(b^{-1/2})$.

Next define the EE $\Psi_n(x)$ for $T_{2N}$ by its Fourier transform,

$$\hat{\Psi}_n(t) \equiv \int_{-\infty}^{\infty}\exp(\iota tx)\,d\Psi_n(x)$$

(4.27)
$$= \exp(-t^2/2)\Bigg[1 + E(\iota tZ_n)^3/6$$
$$+ b^{-3/2}\sum_{i=a+1}^{b_0}\sum_{j=a+1}^{b_0}(\iota t)EZ_{in}V_{jn}\Gamma_n(i,j)\Bigg],$$

where $b_0 = \lceil N/(2\ell) \rceil$ and for $1 \leq i,j \leq b_0$, $Z_{in} = \ell^{-1/2}\sum_{k=1}^{2\ell}Y_{(i-1)2\ell+k}$, $V_{jn} = -N^{-2}b^2\tilde{S}([(i-1)2\ell+1],[2i\ell \wedge N])/[2\sigma_n^2]$ and $\Gamma_n(i,j) = \sum_{r=1}^{2}(r!)^{-1}\times$

$[\iota t G_n(i,j)]^r$, with $G_n(i,j)$ denoting the sum over all $\{Z_{kn}b^{-1/2}: |i-k| \leq 1, |j-k| \leq 1, 1 \leq k \leq b_0\}$. Also, let $Q_N = -Z_N^{(2)} s_n^{(2)}/[2\sigma_n^2]$, the quadratic part of $T_{2N}$. Using arguments in the proofs of Lemma 4.5 above and Lemma 3.5 of Lahiri [19], one can show that

$$
\begin{aligned}
\left| EB_N(t) \exp\left( \iota t \sum_{j=1}^{a} Z_{jn}/\sqrt{b} \right) \right| \\
\leq C_1 [\exp(-C_2 a/m_3) \mathbb{1}(|t| > \kappa b^{1/2}) \\
+ \exp(-C_3 t^2 a^2 /[bm_3]) \mathbb{1}(|t| \leq \kappa b^{1/2})]
\end{aligned}
\tag{4.28}
$$

for some constants $C_i \equiv C_i(\sigma_\infty^2, \kappa)$, where $B_N(t) = \prod_{j=a+1}^{b_0} \{\prod_{k=1}^{r_j}(1+\alpha_{kj}Z_{jn}+\beta_{kj}V_{jn})\} \exp(\iota t Z_{jn})$ for some $\alpha_{kj}, \beta_{kj} \in [-1,1]$ and $r_j \in \mathbb{Z}_+$ with $\sum_{j=a+1}^{b_0} \sum_{k=1}^{r_j} [|\alpha_{kj}| + |\beta_{kj}|] \leq 4$. By Taylor's expansion, for $t \in \mathbb{R}$,

$$
\begin{aligned}
E\exp(\iota t T_{2N}) &= E\exp(\iota t Z_N)\bigg[1 + \iota t Q_N \\
&\quad + 2^{-1}(\iota t Q_N)^2 \int_0^1 \exp(\iota u t Q_N)\, du \bigg] \\
&= \hat{\Psi}_n(t) + R_{3n}(t),
\end{aligned}
\tag{4.29}
$$

where, using (4.28) and Lemma (3.30) of [8], one can show that $|R_{3n}(t)| \leq C(|t|^2 + |t|^4) b^{-1/2-\delta}$ for some $\delta \in (0, 1/2)$.

Now using Esseen's smoothing inequality (cf. Chapter 15, Feller [5]) over the interval $\{|t| < b^{1/2} \log n\}$ and using (4.29) for $|t| \leq (\log n)^4$, the second term on the right-hand side of (4.28) for $|t| \in ((\log n)^4, \kappa b^{1/2}]$, and the first term on the right-hand side of (4.28) for $|t| > \kappa b^{1/2}$, one can conclude that

$$
\sup_{x \in \mathbb{R}} |P(T_N \leq x) - \Psi_n(x)| = o(b^{-1/2}).
\tag{4.30}
$$

Next we derive an EE for the bootstrapped Studentized statistic. Let $\tilde{Z}_{in}^* = (\ell_1/\ell)^{1/2}[\bar{Y}_{in,1}^* - \hat{\mu}_{n,1}]$, $Z_{in}^* = \tilde{Z}_{in}^*/\check{\sigma}_n$ and $Z_N^* = b_1^{-1/2} \sum_{i=1}^{b_1} Z_{in}^*$, where $\bar{Y}_{in,1}^*$ is the average of the $i$th resampled BOBB block (of size $\ell_1$), $\hat{\mu}_{n,1} = E_*\bar{Y}_{1n,1}^*$ and $\check{\sigma}_n^2 = \mathrm{Var}_*(\tilde{Z}_{in}^*)$. Using Taylor's expansion and moderate deviation inequalities for sums of independent random variables, one can show that

$$
\begin{aligned}
T_N^* &= Z_N^* - Z_N^*\bigg[b_1^{-1}\sum_{j=1}^{b_1}\{(\tilde{Z}_{jn}^*)^2 - \check{\sigma}_n^2\}\bigg]\bigg/(2\check{\sigma}_n^2) + R_{1N}^* \\
&\equiv T_{1N}^* + R_{1N}^*, \qquad \text{say,}
\end{aligned}
$$

where $A_n \equiv \{\check{\sigma}_n^2 > \sigma_\infty^2/2\} \cap \{E_*|Z_N^*|^{r_0} < C\}$ and $P_*(|R_{1N}^*| > b^{-1/2}(\log n)^{-1}) \leq Cb^{-1/2}(\log \log n)^{-1}$. Next define $T_{2N}^*$ by deleting the first $a_1 \equiv \lceil b_1^{(1-2\alpha)} \rceil$-many $Z_{in}^*$'s from the quadratic term in $T_{1N}^*$, where $\alpha = \frac{1}{2} - \frac{(1-5\delta_0)}{6}$. Let $R_{2N}^* = T_{1N}^* - T_{2N}^*$. Then, with $\varepsilon_n = \frac{\sqrt{b_1}}{\sqrt{b}\log n}$, using arguments similar to the unbootstrapped case, one can show that on the set $A_n$,

$$(4.31) \qquad P_*(|R_{2N}^*| > 3\varepsilon_n b_1^{-1/2}) \leq 3a_1\{\varepsilon_n^2 b_1\}^{-1}.$$

Here the choices of $a_1$ and $\varepsilon_n$ ensure that $[\varepsilon_n b_1^{-1/2}] \vee [\frac{a_1}{\varepsilon_n^2 b_1}] = O(b^{-1/2}(\log n)^{-1})$. Next define the EE $\check{\Psi}_n(x)$ for $T_{2N}^*$ (and also for $T_N^*$) by its Fourier transform,

$$\begin{aligned}
\hat{\check{\Psi}}_n^*(t) &\equiv \int_{-\infty}^{\infty} \exp(\iota t x)\, d\check{\Psi}_n(x) \\
&= \exp(-t^2/2)\bigg[1 + E_*(\iota t Z_{1n}^*)^3/[6\sqrt{b_1}] \\
&\qquad + b_1^{-3/2} \sum_{i=a_1+1}^{b_1} \sum_{j=a_1+1}^{b_1} (\iota t) E(\tilde{Z}_{in}^*[(\tilde{Z}_{in}^*)^2 - \check{\sigma}_n^2] \\
&\qquad\qquad \times \{-2\check{\sigma}_n^3\}^{-1} \Gamma_n^*(i,j))\bigg],
\end{aligned}$$
(4.32)

where $\Gamma_n^*(i,j) = \sum_{r=1}^{2}(r!)^{-1}[\iota t G_n^*(i,j)]^r$, with $G_n(i,j)$ denoting the sum over all $\{Z_{kn}^* b^{-1/2} : k \in \{i,j\}\}$. Using the independence of the $Z_{in}^*$'s, the moment condition and Taylor's expansion, one can show (cf. Götze [6]) that on the set $A_n$,

$$(4.33) \qquad \int_{|t| \leq m_3} \frac{|E_* \exp(\iota t T_{2N}^*) - \hat{\check{\Psi}}_n^*(t)|}{|t|}\, dt \leq C b_1^{-1} \int_{|t| \leq m_3} \frac{[t^2 + t^4]}{|t|}\, dt \\
\leq C b_1^{-1} m_3^4.$$

Since the integral of $|\hat{\check{\Psi}}_n^*(t)||t|^{-1}$ over $\{|t| > m_3\}$ is $O(\exp(-Cm_3^2))$, in view of Esseen's lemma, it is enough to consider the integral of $|E_* \exp(\iota t T_{2N}^*)||t|^{-1}$ over $\{t: m_3 < |t| < b_1\}$ (say). Let $Q_N^*$ denote the quadratic part in $T_{2N}^*$. Then

$$\begin{aligned}
&|E_* \exp(\iota t T_{2N}^*)| \\
&= E_*\bigg[\exp(\iota t Z_N^*)\bigg\{\sum_{r=0}^{r_0-1} \frac{[\iota t Q_N^*]^r}{r!} \\
&\qquad + \frac{[\iota t Q_N^*]^{r_0}}{r_0!} \int_0^1 E_* \exp(\iota u t Q_N^*)\, du\bigg\}\bigg].
\end{aligned}$$
(4.34)

Note that by the independence of the $\tilde{Z}_{kn}^*$'s, the $r$th summand above for $r \in \{0, \ldots, r_0 - 1\}$ is bounded above by $C_1(r)(b_1)^r |t|^r \exp(-C_2(r)t^2)$, while the last term is bounded above by $C_1 |t|^{r_0} \exp(-C_2 a_1 t^2 / b_1)$. Hence

$$(4.35) \quad \int_{m_3 < |t| \leq b_1^\alpha m_3} |E_* \exp(\iota t T_{2N}^*)| \frac{dt}{|t|} \leq C[b_1]^{-(1/2-\alpha)r_0}[m_3]^{r_0} = o(b^{-1/2}).$$

Next, using the independence of $\{Z_{kn}^* : k = 1, \ldots, a_1\}$ and $Q_N^*$, one can show that $|E_* \exp(\iota t T_{2N}^*)| \leq \exp(-Ct^2 a_1 / b_1)$ for all $|t| \leq \kappa b_1^{1/2}$, for some $\kappa > 0$, so that

$$(4.36) \quad \begin{aligned} \int_{b_1^\alpha m_3 < |t| < \kappa b_1^{1/2}} &|t|^{-1} |E_* \exp(\iota t T_{2N}^*)| \, dt \\ &\leq C_1 \frac{[b_1]^{1-2\alpha}}{a_1 [m_3]^2} \exp(-C_2 m_3^2 a_1 [b_1]^{2\alpha-1}) \leq C_3 m_3^{-2} \exp(-C_4 m_3^2). \end{aligned}$$

Finally, consider $|t| \in [\kappa b_1^{1/2}, b^{1/2} \log n]$. For any $\varepsilon > 0$,

$$\sup\{|E \exp(\iota t \tilde{Z}_{in}^*)| : \varepsilon < |t| < b_1\}$$

$$(4.37) \quad \begin{aligned} &\leq \sup\left\{\left|N_1^{-1} \sum_{i=1}^{N_1} [\exp(\iota t \bar{Y}_{in,1}(\ell_1/\ell)^{1/2}) \right.\right. \\ &\quad \left.\left. - E \exp(\iota t \bar{Y}_{in,1}(\ell_1/\ell)^{1/2})]\right| : \varepsilon < |t| < b_1\right\} \\ &\quad + \sup\{|E \exp(\iota t \bar{Y}_{1n,1}(\ell_1/\ell)^{1/2})| : \varepsilon < |t| < b_1\} \\ &\equiv I_{1N} + I_{2N}, \quad \text{say.} \end{aligned}$$

Using a discretizing argument as in Lemma 4.2 of Babu and Singh [1], one can show (cf. Lahiri [19]) that

$$(4.38) \quad P(I_{1N} > (\log[n+1])^{-1}) \to 0 \quad \text{as } n \to \infty.$$

Note that $\bar{Y}_{1n,1}(\ell_1/\ell)^{1/2} = (\ell_1/\ell)^{-1/2} \sum_{i=1}^{\ell_1/\ell} \bar{Y}_{in}$, where $\bar{Y}_{in} = \ell^{-1} \sum_{j=i}^{i+\ell-1} Y_{jn}$ [cf. (2.1)]. Next, using Theorem 4.1 [with $f(\cdot) = \exp(\iota t(\cdot))$] for $\kappa \leq |t| \leq (\ell_1/\ell)^{1/2}$ and Lemma 4.5 and condition (iv) for $|t| \geq (\ell_1/\ell)^{1/2}$, one can show that there exists a $\kappa \in (0, 1)$ such that

$$(4.39) \quad I_{2N} \leq (1 - \kappa)$$

for all $n > \kappa^{-1}$. Hence, by (4.31)–(4.39), it follows that

$$(4.40) \quad \sup_{x \in \mathbb{R}} |P_*(T_N^* \leq x) - \check{\Psi}_n(x)| = o_P(b^{-1/2}).$$

Next comparing the EEs for $T_N$ and $T_N^*$, using the blocking arguments as in derivation of (4.25) and observing that $b_1^{-1/2} E\{(\ell_1/\ell)^{1/2} \bar{Y}_{1n,1}\}^3 = [b_1^{-1/2} \ell_1/\ell)^{-1/2}] \times [(\ell_1/\ell)^{-1} E\{\sum_{i=1}^{\ell_1/\ell} Y_{in}\}^3] = b^{-1/2}[b^{-1} E\{\sum_{i=1}^{b} Y_{in}\}^3](1+o(1))$, one can complete the proof of (3.8). □

DEPARTMENT OF STATISTICS
IOWA STATE UNIVERSITY
AMES, IOWA 50011
USA
E-MAIL: snlahiri@iastate.edu